\documentclass{article}

\usepackage{amssymb}
\usepackage{latexsym}
\usepackage{theorem}

\newcommand{\be}{\begin{equation}}
\newcommand{\ee}{\end{equation}} 
\newcommand{\bea}{\begin{eqnarray}}
\newcommand{\beann}{\begin{eqnarray*}}
\newcommand{\eea}{\end{eqnarray}}
\newcommand{\eeann}{\end{eqnarray*}}

\newcommand{\la}{\langle}
\newcommand{\ra}{\rangle}

\newcommand{\pa}{\partial}
\newcommand{\al}{\alpha}
\newcommand{\bt}{\beta}

\makeatletter
\@addtoreset{equation}{section}
\makeatother

\theoremstyle{break}
\newtheorem{th1}{Theorem}[section]
\newtheorem{p1}{Proposition}[section]

\newtheorem{thth1}[th1]{Theorem}
\newtheorem{thp1}[p1]{Theorem}

\newtheorem{pp1}[p1]{Proposition}

\newtheorem{cp1}[p1]{Corollary}

\newtheorem{lp1}[p1]{Lemma}

\begin{document}

\begin{center}
{\Large\bf The fake monster superalgebra}\\[1cm]
Nils R. Scheithauer\footnote{
Diese Arbeit wurde mit Unterst\"utzung eines Stipendiums im Rahmen des 
Gemeinsamen Hochschulsonderprogramms III von Bund und L\"andern \"uber den 
DAAD erm\"oglicht.}\\
D.P.M.M.S.\\
16 Mill Lane\\
Cambridge CB2 1SB\\
England\\
\end{center}
\vspace*{2cm}

\noindent
We show that the physical states of a 10 dimensional superstring moving on a torus form a generalized Kac-Moody superalgebra. This gives the first explicit realizations of these algebras. For a special torus the denominator function of this algebra is an automorphic form so that we can determine the simple roots. We call this algebra the fake monster superalgebra. 

\section{Introduction}

The physical states of a bosonic string moving on a torus form a nice Lie algebra called the fake monster algebra (\cite{B2}). The study of this algebra has led to the definition of vertex algebras, generalized Kac-Moody algebras and the proof of the moonshine conjectures (\cite{B4}). The fake monster algebra is the simplest example of a generalized Kac-Moody algebra. It is easy to extend the definition of generalized Kac-Moody algebras to superalgebras. Since the physical states of a compactified bosonic string form a generalized Kac-Moody algebra one would expect that the physical states of a compactified superstring give a realization of a generalized Kac-Moody superalgebra. Some evidence for the existence of such an algebra is given in \cite{NRS} and \cite{B4}. In this paper we show that this conjecture is true. In a special case we get a particularly nice algebra which we call the fake monster superalgebra.

We explain the construction in more detail. The superstring has bosonic and fermionic degrees of freedom. The bosonic degrees are represented by a vertex algebra constructed from a lattice. The fermionic degrees can be bosonized and then also be described by such an algebra. One also has to include some ghost coordinates to define a fermion emission vertex. Thus the Fock space of a compactified superstring can be described by the vertex algebra of some suitable lattice. This vertex algebra carries representations of various algebras the most important being the $N=2$ superconformal algebra. The representation of this algebra is used to construct the BRST operator $Q$ satisfying $Q^2=0$. The physical states of the superstring are now represented by the cohomology of this operator. More precisely for nonzero momentum there are infinitely many isomorphic copies of vector spaces called pictures representing the physical states. We restrict to the canonical ghost pictures. The vertex algebra induces a product on the vector space of physical states that defines the structure of a Lie superalgebra on this space. This Lie superalgebra has a number of nice properties. It is graded by the 10 dimensional momentum lattice. The dimensions of the root spaces can be described easily. All roots have zero or negative norm. It has an invariant supersymmetric bilinear form and carries a representation of the $N=1$ supersymmetry algebra. In particular it is a generalized Kac-Moody superalgebra.
 If the momentum lattice of the superstring is the unique unimodular even 10 dimensional Lorentzian lattice then the denominator function is an automorphic form with known product expansion. This allows us to determine the simple roots in this case.

Most of the ideas necessary for the above construction come from physics and are described in the physics literature (cf. \cite{FMS},\cite{LT},\cite{K} and \cite{LZ2}). References \cite{B5} and \cite{R2} are also important for this work.

We describe the sections of this paper in more detail.

In the first section we recall some facts about vertex algebras graded by abelian groups. We show that in the lattice construction of these algebras the locality axiom is a simple consequence of Wick's Theorem.

In the second section we describe some results about generalized Kac-Moody superalgebras and state the characterization theorem that will be applied to deduce the first main result of this paper.

In the next section we construct the vertex algebra of the 10 dimensional superstring moving on a torus from a rational 18 dimensional lattice. We give the operator product expansions of various fields. They are used to construct an action of the $N=2$ superconformal algebra and to construct the nilpotent BRST operator.

In the last section we describe the cohomology spaces that represent the physical states and define a Lie bracket on the space of physical states. We derive a number of properties of this Lie superalgebra, e.g. a root space decomposition, the multipicities of the roots, an action of the $N=1$ supersymmetry algebra etc.. These properties imply that this algebra is a generalized Kac-Moody superalgebra. Finally we define the fake monster superalgebra which is a special case of the above construction with the Lorentzian lattice $II_{9,1}$ as momentum lattice and determine the simple roots of this algebra using the product expansion of its denominator function.

\section{Vertex algebras}

In this section we recall some results on vertex algebras. A more detailed exposition can be found in \cite{K2} and \cite{DL}. We give a simple proof of the locality axiom in the lattice construction.

\subsection{Definition and some properties} 

Let $\Gamma$ be a finite abelian group and $\Delta : \Gamma \times \Gamma \rightarrow \mathbb{Q}/\mathbb Z$ a symmetric map bilinear mod $\mathbb Z$. This implies that $\Delta$ takes it values in $\frac{1}{g}\mathbb Z$ where $g$ is an exponent of $\Gamma$. Let $\eta \, :\, \Gamma \times \Gamma \rightarrow \mathbb C^*$ be bimultiplicative. Then $\eta(\gamma_1,\gamma_2)$ is a root of unity, i.e. $\eta(\gamma_1,\gamma_2)^g=1$.

Suppose that 
\be V=\bigoplus_{\gamma \in \Gamma} V_{\gamma} \ee 
is a $\Gamma$-graded 
complex vector space and 
\bea 
V & \rightarrow & (\mbox{End} V)[[z^{\frac{1}{g}},z^{-\frac{1}{g}}]] \\
a & \mapsto     & a(z)=\sum_{n\in \frac{1}{g}\mathbb Z} a_nz^{-n-1}
\eea
a parity preserving state-field correspondence, i.e. for $a\in V_{\gamma_1}, 
b\in V_{\gamma_2}$ 
\be \label{dpe}
a(z)b = \sum_{n\in {\mathbb Z}+\Delta(\gamma_1,\gamma_2)} a_nb\, z^{-n-1} 
\ee
and
\be
a_nb \in  V_{\gamma_1+\gamma_2} 
\ee
is zero for $n$ sufficiently large. 

The field $a(z)$ is also written $Y(a,z)$ in the literature.

$V$ is a $\Gamma$-graded vertex algebra if it satisfies the following 
conditions
\begin{enumerate}
\item There is an element $1\in V_0$, called vacuum, with 
      \be 1(z)a=a \qquad \mbox{and} \qquad a(z)1\big|_{z=0}=a. \ee 
\item The operator $D$ on $V$ defined by $Da=a_{-2}1$ satisfies 
      \be [D,a(z)]=\pa a(z). \ee
\item The locality condition 
      \be
      i_{z,w}(z-w)^n\, a(z)b(w)-\eta(\gamma_1,\gamma_2)i_{w,z}(z-w)^n \, b(w)a(z)=0
      \label{loc} \ee
      holds for $a\in V_{\gamma_1}, b\in V_{\gamma_2}$ and 
      $n\in \mathbb Z+\Delta(\gamma_1,\gamma_2)$ sufficiently large.
      Here $i_{z,w}(z-w)^n$ is the binomial expansion of $(z-w)^n$ in the 
      domain $|z|>|w|$, i.e. 
      $i_{z,w}(z-w)^n=\sum_{k\geq 0}(-1)^k{n \choose k}z^{n-k}w^k$.
      In the expresssion $i_{w,z}(z-w)^n$ we must be precise about the roots of      unity. We choose $i_{w,z}(z-w)^n=e^{i\pi n}i_{w,z}(w-z)^n=
      e^{i\pi n}\sum_{k\geq 0}(-1)^k{n \choose k}w^{n-k}z^k$.
\end{enumerate} 
Note that $D$ leaves the spaces $V_{\gamma}$ invariant.

For the rest of the section assume that $V$ is a $\Gamma$-graded vertex 
algebra. 

We now describe some consequences of the definition. 

Locality implies $\eta(\gamma_1,\gamma_2)e^{\Delta(\gamma_1,\gamma_2)}=
                  \eta(\gamma_2,\gamma_1)^{-1}e^{-(\gamma_2,\gamma_1)}$.

For all $a$ in $V$ one has
\bea 
a(z)1   &=& e^{zD}a \\
(Da)(z) &=& \pa \, a(z) \label{Dpa}\\
e^{wD} a(z) e^{-wD} &=& a(z+w) \quad\mbox{in the domain }|z|>|w|.
\eea

For $a\in V_{\gamma_1}$ and $b\in V_{\gamma_2}$ the following symmetry formula can be proved as in the superalgebra case.
\be 
a(z)b=\eta(\gamma_1,\gamma_2)e^{zD}b(-z)a
\, .
\ee

The most important property of vertex operators is Borcherds' identity.
\begin{p1}
Let $a\in V_{\gamma_1},\, b\in V_{\gamma_2},\, c\in V_{\gamma_3}$ and 
$n\in \mathbb Z+\Delta(\gamma_1,\gamma_2),\, k\in {\mathbb Z}+\Delta(\gamma_1,\gamma_3)$. Then we have 
\[ \sum_{j\geq 0} {n \choose j} (-1)^j 
          \big\{ a_{n+k-j}b_{m+j}c  
         -\eta({\gamma_1},{\gamma_2})e^{i\pi n} b_{m+n-j}a_{k+j}c \big\} \]  
\[ = \sum_{j\geq 0} {k \choose j} (a_{n+j}b)_{k+m-j}c \]
\end{p1}

We will use some special cases of Borcherds' identity very often in the following sections. If $\Delta(\gamma_1,\gamma_3)$ is zero we can put $k=0$ to get the associativity formula 
\be 
(a_{n}b)_{m}c =
\sum_{j\geq 0} {n \choose j} (-1)^j 
          \big\{ a_{n-j}b_{m+j}c  
         -\eta({\gamma_1},{\gamma_2})e^{i\pi n} b_{m+n-j}a_{j}c \big\}.
\ee
If $\Delta(\gamma_1,\gamma_2)=0$ the commutator formula 
\be  
a_{k}b_{m}c  
         -\eta({\gamma_1},{\gamma_2})b_{m}a_{k}c =
\sum_{j\geq 0} {k \choose j} (a_{j}b)_{k+m-j}c 
\ee 
holds.
Now suppose that $a_{n+j}b=0$ for $j\geq 1$. Then 
\be \label{bisc}
\! \! (a_{n}b)_{k+m}c = 
\sum_{j\geq 0} {n \choose j} (-1)^j 
          \big\{ a_{n+k-j}b_{m+j}c  
         -\eta({\gamma_1},{\gamma_2})e^{i\pi n} b_{m+n-j}a_{k+j}c \big\}. 
\ee
In contrast to the associativity formula we need not assume here that $\Delta(\gamma_1,\gamma_3)$ is zero.

We conclude this section with some remarks.
\pagebreak[3]

An element $\omega\in V_{0}$ is called Virasoro element of central charge 
 $c$ if it satisfies the following 3 conditions
\begin{enumerate}
\item[] The operators $L_m=\omega_{m+1}$ give a representation of the 
          Virasoro algebra of central charge $c$, i.e.
          \[ [L_m,L_n]=(m-n)L_{m+n}+ \frac{m^3-m}{12}\delta_{m+n,0}\,c\, . \] 
\item[] $L_0$ is diagonizable on $V$.
\item[] $D=L_{-1}$.
\end{enumerate}

In conformal field theory the vertex operator of a state $a$ with $L_0 a=h a$  is expanded as 
$a(z)=\sum  a_{(n)}z^{-n-h}$, i.e. the expansion of a field depends on its conformal weight $h$. The relation between the two modes is given by 
$a_n= a_{(n+1-h)}$ resp. $a_{(n)}=a_{n-1+h}$. 

\subsection{Construction from rational lattices} \label{lc}

Now we describe how vertex algebras can be constructed from rational lattices.
The proof given here only uses Wick's theorem (cf. \cite{K2}).
 
Let $L$ be a rational lattice of finite rank and $L_0\neq 0$ an even sublattice of $L$ with $L\subset L_0^*$. Then $(L,L_0)\subset \mathbb Z\,$.
Let $L=L_0\cup L_1\cup \ldots \cup L_n$ with $L_i=\delta_i + L_0$ be the coset decomposition of $L$ with respect to $L_0$. Put $\Gamma=L/L_0=\{\gamma_0,\gamma_1,\ldots , \gamma_n\}$ where $\gamma_i$ corresponds to $L_i$ and let $g$ be an exponent of $\Gamma$. 
Define $\Delta\, :\, \Gamma \times \Gamma \rightarrow \mathbb Q/\mathbb Z$ by 
$\Delta(\gamma_i,\gamma_j)=-(\delta_i,\delta_j) \, \mbox{mod}\, \mathbb Z$.
This map is well-defined and bilinear mod $\mathbb Z\,$. 
Let $\eta \, :\, \Gamma \times \Gamma \rightarrow \mathbb C^*$ be a bimultiplicative map with
\be \eta(\gamma_j,\gamma_j)=e^{i\pi (\delta_j,\delta_j)}\, . \label{nc} \ee
This condtion is quite restrictive but necessary for consistency.

Let $\varepsilon : L\times L\rightarrow {\mathbb C}^*$ be a 2-cocycle satisfying $\varepsilon(\alpha,\beta)=B(\alpha,\beta)\varepsilon(\beta,\alpha)$ where
$B(\alpha,\beta)=e^{-i\pi (\alpha,\beta)}\eta(\gamma_i,\gamma_j)$ for $\alpha \in L_i$ and $\beta \in L_j$ and $\varepsilon(0,\alpha)=\varepsilon(\alpha,0)=1$. Note that $B$ is a bimultiplicative map from $L\times L$ to ${\mathbb C}^*$ and that $B(\alpha,\alpha)=1$ by (\ref{nc}).

It is easy to construct such a 2-cocycle. Let $\{\alpha_1,\ldots, \alpha_m\}$ be a $\mathbb Z$-basis of $L$. Define 
\[ \begin{array}{lclr} 
\varepsilon(\alpha_i,\alpha_j) &=& a_{ij} & i\leq j \\
\varepsilon(\alpha_i,\alpha_j) &=& B(\alpha_i,\alpha_j)
    \varepsilon(\alpha_j,\alpha_i)   & i>j 
\end{array} \]
with some $a_{ij}\in {\mathbb C}^*$ and extend to $L$ by bimulticativity. This gives us a 2-cocycle $\varepsilon : L\times L\rightarrow {\mathbb C}^*$ satisfying $\varepsilon(\alpha,\beta)=B(\alpha,\beta)\varepsilon(\beta,\alpha)$ and $\varepsilon(0,\alpha)=\varepsilon(\alpha,0)=1$.

Extend the symmetric form of $L$ to $h=L{\otimes}_{\mathbb Z}{\mathbb C}$. Define the infinite dimensional Heisenberg algebra $\hat{h}=h\otimes {\mathbb C}[t,t^{-1}]\oplus {\mathbb C}c$ with products $[h_1(m),h_2(n)]=m\delta_{m+n,0} \la h_1,h_2 \ra c$ and $[h_1(m),c]=0$. Then $\hat{h}^-=h\otimes t^{-1}{\mathbb C}[t^{-1}]$ is an abelian subalgebra of $\hat{h}$ and $S(\hat{h}^-)$ is the 
symmetric algebra of polynomials in $\hat{h}^-$.  
  
Let ${\mathbb C}[L]$ be the group algebra of $L$ with basis $\{ e^{\alpha} |\, \alpha \in L \}$ and products $e^{\alpha}e^{\beta}=e^{\alpha + \beta}$.

The vector space 
\be V=S(\hat{h}^-)\otimes {\mathbb C}[L] \ee
decomposes as
\be V=\bigoplus_{\gamma_i \in \Gamma} V_{\gamma_i}  \ee
where 
\be V_{\gamma_i}=S(\hat{h}^-)\otimes {\mathbb C}[L_i]\, . \ee

There is a natural action of $\hat{h}$ on $V$. 

We define the vertex operator of $e^{\alpha}$ as
\be 
e^{\alpha}(z) = e^{\alpha}(z)^{+}\, e^{\alpha}(z)^{-}
\ee
where
\be e^{\alpha}(z)^{+}=e^{\al}c_{\alpha}e^{\sum_{m>0}{\alpha}(-m)\frac{z^m}{m}}
                     =e^{\al}c_{\alpha}\sum_{m\geq 0}S_m(\al)z^m \ee
and
\be
  e^{\alpha}(z)^{-}=z^{\al(0)}e^{-\sum_{m>0}{\alpha}(m)\frac{z^{-m}}{m}} \, .
\ee
The linear operator $c_{\alpha}$ acts on ${\mathbb C}[L]$ as $c_{\alpha} e^{\beta}=\varepsilon(\alpha,\beta) e^{\beta}$. 

For $h(-n-1),\, n\geq 0$ put
\be h(-n-1)(z) = \pa_z^{(n)} h(z) \ee
with $h(z)=\sum_{n\in {\mathbb Z}}h(n)z^{-n-1}$ and 
$\pa_z^{(n)}=\frac{\pa_z}{n!}$.

The bosonic normal ordering $: h_1(n_1)\cdots h_k(n_k) :$ of Heisenberg generators is defined by putting all $h(n)$ with $n<0$ ("creation operators") to the left of those with $n\geq 0$ ("annihilation operators").

Define 
\bea 
\lefteqn{ \left( h_1(-n_1-1)\cdots h_k(-n_k-1) e^{\alpha}\right)(z)  } 
\nonumber \\
&= e^{\alpha}(z)^{+}
    \, : h_1(-n_1-1)(z) \cdots h_k(-n_k-1)(z) : \, 
   e^{\alpha}(z)^{-} &
\eea
and extend this definition linearly to $V$ to get a parity preserving state-field correspondence $V \rightarrow (\mbox{End} V)[[z^{\frac{1}{g}},z^{-\frac{1}{g}}]]$.

\begin{thp1}
With this structure $V$ is a vertex algebra graded by $\Gamma$. The vacuum is given by $1\otimes e^0$.
\end{thp1}
{\it Proof:}
We only prove the locality. Using
\be
e^{\alpha}(z)^{-}e^{\beta}(w)^{+} =
(z-w)^{(\alpha,\beta)}e^{\beta}(w)^{+}e^{\alpha}(z)^{-} 
                                   \quad |z|>|w| 
\ee
we get for $\al\in L_i$ and $\bt \in L_j$ in the domain $|z|>|w|$
\[ e^{\alpha}(z) e^{\beta}(w) = 
    (z-w)^{(\alpha,\beta)}e^{\alpha}(z)^{+} e^{\beta}(w)^{+}
                          e^{\alpha}(z)^{-} e^{\beta}(w)^{-} \, . \]
While the commutator $[e^{\alpha}(z)^{-},e^{\beta}(w)^{-}]$ is zero 
$e^{\alpha}(z)^{+}$ and $e^{\beta}(w)^{+}$ do not commute. 

Let $n\in {\mathbb Z}+{\Delta}(\gamma_i,\gamma_j)$. Then
\[ i_{z,w}(z-w)^n e^{\al}(z) e^{\bt}(w)
          -\eta(\gamma_i,\gamma_j)i_{w,z}(z-w)^n e^{\bt}(w) e^{\al}(z)=0 \]
is equivalent to 
\[ e^{\alpha}c_{\al}e^{\bt}c_{\bt}i_{z,w}(z-w)^{n+(\alpha,\beta)}
   -\eta(\gamma_i,\gamma_j)e^{-i\pi (\al,\bt)}
   e^{\bt}c_{\bt}e^{\al}c_{\al}i_{w,z}(z-w)^{n+(\alpha,\beta)}=0 \]

Note that $n+(\alpha,\beta)$ is an integer by the definition of $\Delta$. For $n$ sufficiently large this integer is positive. Then 
$i_{z,w}(z-w)^{n+(\alpha,\beta)}=i_{w,z}(z-w)^{n+(\alpha,\beta)}$. Hence 
$e^{\alpha}(z)$ and $e^{\beta}(w)$ are local if 
\be 
e^{\alpha}c_{\al}e^{\bt}c_{\bt}=B(\al,\bt)e^{\bt}c_{\bt}e^{\alpha}c_{\al} \, .
\ee
Since $B$ is bimultiplicative this equation is equivalent to the equations
\bea
\varepsilon(\al,\bt) &=& B(\al,\bt) \varepsilon(\bt,\al) \\
\varepsilon(\al,\bt)\varepsilon(\al+\bt,\gamma) &=&
\varepsilon(\bt,\gamma)\varepsilon(\al,\bt+\gamma) \, . 
\eea
But these are satisfied by construction of $\varepsilon$ and $B$. It follows that the fields $e^{\alpha}(z)$ and $e^{\beta}(w)$ are local.

Now we want to prove the locality of the fields of two arbitrary states in $V$.
For that it is sufficient to prove it for two states of the form
\beann
a &=& h_1(-n_1-1)\cdots h_i(-n_i-1) e^{\al} \\
b &=& k_1(-m_1-1)\cdots k_j(-m_j-1) e^{\bt}\, . 
\eeann
We do this by induction on $i+j$. The case of $i+j=0$ has been done above.
We can assume that $j\geq 1$ so that we can write 
$b=c\otimes d$ with $c=k_1(-m_1-1)$.
Then 
\[ b(w)=c(w)^{+} d(w) + d(w) c(w)^{-} \]
where 
\[ c(w)^{-}=k_1(-m_1-1)(w)^{-}
=\sum_{l\geq 0}(-1)^l {l \choose m_1}k_1(l-m_1)w^{-l-1} \]
only contains annihilation operators and 
\[ c(w)^{+}=k_1(-m_1-1)(w)^{+}
=\sum_{l< 0}(-1)^l {l \choose m_1}k_1(l-m_1)w^{-l-1}\]
only creation operators. We get
\bea
\lefteqn{ 
i_{z,w}(z-w)^n a(z) b(w) -\eta(\gamma_i,\gamma_j)i_{w,z}(z-w)^nb(w)a(z)} 
\label{ind} \\
&=& c(w)^{+}\big\{i_{z,w}(z-w)^n a(z) d(w)
     -\eta i_{w,z}(z-w)^n d(w)a(z) \big\} \nonumber \\
& & +\big\{i_{z,w}(z-w)^n a(z) d(w)
     -\eta i_{w,z}(z-w)^n d(w)a(z) \big\}
     c(w)^{-} \nonumber \\
& & +\big\{i_{z,w}(z-w)^n [a(z),c(w)^{+}]d(w)
     -\eta i_{w,z}(z-w)^n d(w)[c(w)^{-},a(z)]\big\} \, . \nonumber
\eea
The first two terms vanish for $n\in {\mathbb Z}+{\Delta}(\gamma_i,\gamma_j)$
sufficiently large by the induction hypothesis. Hence we need to calculate 
the commutators. To simplify notations we write $a_l=h_l(-n_l-1)$. Then
$a(z)=e^{\alpha}(z)^{+}:a_1(z)\cdots a_i(z): e^{\alpha}(z)^{-}$ and
\bea
\lefteqn{ [a(z),c(w)^{+}] } \\
&=& e^{\alpha}(z)^{+}:a_1(z)\cdots a_i(z):[\, e^{\alpha}(z)^{-},c(w)^{+}]
    \nonumber \\
& & +\, e^{\alpha}(z)^{+}[\,:a_1(z)\cdots a_i(z):\, ,c(w)^{+}]\,  
    e^{\alpha}(z)^{-} \, . \nonumber
\eea
Using $[e^{\alpha}(z)^{-},h(l)]=-z^l(\al,h)e^{\alpha}(z)^{-}$ for $l<0$ 
we get
\[ [e^{\alpha}(z)^{-},c(w)^{+}]=-i_{z,w}(z-w)^{-m_1-1}(\al,k_1)e^{\alpha}(z)^{-}
       \, . \]
Wick's theorem (cf. \cite{K2}) yields 
\[ 
[:a_1(z)\cdots a_i(z):\, ,c(w)^{+}]=\sum_{l=1}^{i}[a_l(z)^{-},c(w)]\, 
\hat{a}_l(z) \]
where $\hat{a}_l(z)$ is $:a_1(z) \cdots a_i(z):$ with $a_l(z)$ dropped.
The commutator relations of the Heisenberg algebra imply
\[ [a_l(z)^{-},c(w)]=(-1)^{n_l }{ n_l+m_1+1 \choose m_1} (n_l+1) (h_l,k_1)
    i_{z,w}(z-w)^{-n_l-m_1-2} \]
so that
\bea
\lefteqn{ [a(z),c(w)^{+}] } \label{com1}\\
&=& -i_{z,w}(z-w)^{-m_1-1}(\al,k_1)a(z) \nonumber \\
& & +\sum_{l=1}^{i}(-1)^{n_l }{ n_l+m_1+1 \choose m_1} (n_l+1) (h_l,k_1)
    i_{z,w}(z-w)^{-n_l-m_1-2} (\hat{a}_le^{\al})(z) \, .\nonumber
\eea
Note that in the case $i=0$ the second term is zero. 
$[c(w)^{-},a(z)]$ can be calculated similarly. We find
\bea
\lefteqn{[c(w)^{-},a(z)] } \label{com2}\\
&=& -i_{w,z}(z-w)^{-m_1-1}(\al,k_1)a(z) \nonumber \\
& & +\sum_{l=1}^{i}(-1)^{n_l }{ n_l+m_1+1 \choose m_1} (n_l+1) (h_l,k_1)
    i_{w,z}(z-w)^{-n_l-m_1-2} (\hat{a}_le^{\al})(z) \, .\nonumber \,
\eea
Again in the case $i=0$ the second term vanishes.

Putting (\ref{com1}) and (\ref{com2}) into the last term of (\ref{ind})
we see that also this term is zero for $n$ sufficiently large by the induction hypothesis. This finishes the proof of the theorem. \hspace*{\fill} $\Box$\\

Note that $e^{\al}(z)^{+}=\sum_{n\geq 0}D^{(n)}e^{\al}c_{\al}z^n$.

The following simple formulas are useful
\beann 
k(-1)_{n} &=& k(n) \\ 
e^{\al}_n e^{\bt} &=& 
    \varepsilon(\al,\bt) S_{-n-1-(\al,\bt)}(\al)e^{\al+\bt} \, .
\eeann

The first 3 Schur polynomials are
\beann 
S_0(\al) &=& 1 \\
S_1(\al) &=& \al(-1) \\
S_2(\al) &=& \textstyle{\frac{1}{2}} \left(\al(-1)^2+\al(-2)\right) \, .
\eeann

$V$ contains a one-parameter family of Virasoro elements. 

If $L$ is nondegenerate then $V$ is simple.

\subsection{Invariant bilinear forms} \label{ibf}

We will need some results on bilinear forms on vertex superalgebras. They have been studied in detail in \cite{NRS}. There the even and the odd part of a vertex superalgebra are distinguished by their $L_0$-eigenvalues. Here we will consider a slighty different situation.

Let $V$ be a vertex superalgebra with Virasoro element. Suppose that $L_1$ acts locally nilpotent on $V$ and that the eigenvalues of $L_0$ are all integral. 

Define the adjoint vertex operator of $a\in V$ with $L_0 a=ha$ by  
\be 
a(z)^*= \sum_{n\in \mathbb Z} a_n^* z^{-n-1}
\ee
where
\be
a_n^* = (-1)^h \sum_{m \geq 0} \left( \frac{L_1^m}{m!} a \right)_{2h-n-m-2} 
\, .
\ee
For example $L_n^* = \omega_{n+1}^* = \omega_{1-n} = L_{-n}$. 
The adjoint vertex operators satisfy a number of identities. 

A bilinear form $(\, , \, )$ on $V$ is called invariant if
\be (a_n b,c) = (-1)^{|a||b|} (b, a_n^*c) \, . \ee
An invariant bilinear form is supersymmetric, i.e. 
\be (a,b) = (-1)^{|a||b|} (b,a) \ee
and states with different $L_0$-eigenvalues are orthogonal.

Let $V_h$ be the subspace of states with $L_0$-eigenvalue $h$.
We can define a linear functional $f$ on $V_0$ by putting $f(a)=(1,a)$.
This functional vanishes on $L_1V_1$.
The correspondence between invariant bilinear foms on $V$ and linear functionals on $V_0/L_1V_1$ is injective. It is even surjective. Given an element $f$ of the dual of $V_0/L_1V_1$ then the bilinear form defined by 
$(a,b)=(f \circ \pi)(a_{-1}^*b)$, where $\pi$ is the projection from $V$ to 
$V_0/L_1V_1$, is invariant and induces $f$. This implies 
\begin{thp1}
The space of invariant bilinear forms on $V$ is naturally
isomorphic to the dual of $V_0/L_1V_1$.
\end{thp1}

\section{Generalized Kac-Moody superalgebras}

In this section we describe some results about generalized Kac-Moody superalgebras (cf. \cite{R1}, \cite{R2}, \cite{GN} and also \cite{B1},\ldots,\cite{B4}).
Unless otherwise specified the results will hold over the real or complex numbers.

\subsection{Definitions and properties}

Let $a=(a_{ij})_{i,j \in I}$ be a real square matrix, where $I$ is some finite or countable index set and $\tau$ a subset of $I$, with the following properties
\[ \renewcommand{\arraystretch}{1.3} \begin{array}{ll}
 a_{ij} =a_{ji} \\
 a_{ij} \leq 0                    & \mbox{if} \quad i\neq j \\
 2a_{ij}/a_{ii}\in \mathbb Z      & \mbox{if} \quad a_{ii} >0 \\
  a_{ij}/a_{ii}\in \mathbb Z      & \mbox{if} \quad a_{ii} >0 \; 
                            \mbox{and} \; i\in \tau \, .
\end{array}\]
Let $\hat{G}$ be the Lie superalgebra with generators $e_i, f_i, h_{ij}$ and relations
\[ \renewcommand{\arraystretch}{1.3} \begin{array}{l}
 {[} e_i, f_j ] = h_{ij}  \\
 {[} h_{ij}, e_k ] =  \delta_{ij} a_{ik}e_k  \\
 {[} h_{ij}, f_k ] = -\delta_{ij} a_{ik}f_k  \\
 (\mbox{ad}\, e_i)^{1-2a_{ij}/a_{ii}}e_j = 
 (\mbox{ad}\, f_i)^{1-2a_{ij}/a_{ii}}f_j = 0 
 \quad \mbox{if} \quad a_{ii}>0 \\
 {[} e_i,f_j ] = [f_i,f_j] = 0 \quad \mbox{if} \quad a_{ij} = 0 \\
\mbox{deg } e_i =\mbox{deg } f_i = 0 \quad \mbox{if} \quad i\not\in\tau \\
\mbox{deg } e_i =\mbox{deg } f_i = 1 \quad \mbox{if} \quad i\in \tau 
\end{array} \] 
We list some properties of $\hat{G}$. The elements $h_{ij}$ are zero unless the $i$-th and $j$-th columns of $a$ are equal. The nonzero $h_{ij}$ are linearly independent and span a commutative subalgebra $\hat{H}$ of $\hat{G}$. The elements $h_{ii}$ are usually denoted $h_i$. Every nonzero ideal of $\hat{G}$ has nonzero intersection with $\hat{H}$. The center of $\hat{G}$ is in $\hat{H}$ and contains all the elements $h_{ij}$ with $i\neq j$. 

A Lie superalgebra $G$ is called a generalized Kac-Moody superalgebra if $G$ is the semidirect product $(\hat{G}/C)\cdot D$ where $C$ is a subspace of the center of $\hat{G}$ containing the the nonzero $h_{ij}$ and $D$ an Abelian even subalgebra of $G$ such that the elements $e_i$ and $f_i$ are all eigenvectors of $D$. 

Let $G$ be a generalized Kac-Moody superalgebra. We keep the above notations for the images in the quotient. 

The root lattice $Q$ of $G$ is the free abelian group generated by the elements $\al_i$,$ \, i\in I$, with the bilinear form defined by $(\al_i,\al_j)=a_{ij}$. The elements $\al_i$ are called simple roots. $G$ is graded by the root lattice if we define the degree of $e_i$ as $\al_i$ and the degree of $f_i$ as 
$-\al_i$. We have the usual definitions of roots and root spaces. A root is called real if it has positive norm and imaginary else. 

$H=\hat{H}\oplus D$ is a commutative even subalgebra of $G$ called Cartan subalgebra. There is a unique invariant supersymmetric bilinear form on $G$ satisfying $(h_i,h_j)=a_{ij}$. We have a natural homomorphism of abelian groups from $Q$ to $H$ sending $\al_i$ to $h_i$. By abuse of terminology the images of roots under this map are also called roots. It is important to note that this map is in general not injective. It is possible that $n>1$ imaginary simple roots have the same image $h$ in $H$. In this case we call $h$ a root of multiplicity $n$. If the matrix $a$ has infinitely many identical columns this multiplicity may even be infinite. The root space of root $\al\in Q$ is either even or odd. The root space of the corresponding element in $H$ has in general an even and an odd part. 

Under certain conditions
one can also prove a character formula (\cite{R1},\cite{GN}).

\subsection{A characterization theorem}

While most generalized Kac-Moody algebras posses an almost positive definite contravariant bilinear form in the super case this does not even hold for the finite dimensional algebras. A more general characterization is the following (cf. \cite{R2}, Theorem 3.2) 
\begin{thth1} \label{ut}
Let $G$ be a real Lie superalgebra satisfying the following conditions 
\begin{enumerate}
\item
$G$ has a selfcentralizing even subalgebra $H$, such that $G$ is the direct
sum of the eigenspaces of $H$ and all the eigenspaces are finite dimensional.
The nonzero eigenvalues of $H$ having on $G$ are called roots of $G$.
\item
$G$ has a nonsingular invariant supersymmetric bilinear form.
\item
There is an element $h\in H$ such that the centralizer of $h$ in $G$ is $H$ 
and there exist only finitely many roots $\al$ of $G$ with 
$\left| \al(h) \right| < M$ for any real number $M$
(such an element $h$ is called a regular element). 
A root $\al$ is called positive or negative depending on wether its value on $h$ is positive or negative.
\item
All roots are either of finite type or of infinite type.
(A root $\al$ is of finite type if for all roots $\bt$, $\bt+n\al$ is a root for only finitely many integers $n$. A root $\al$ is of infinite type if, either $\al$ has zero norm and is not of finite type, or for all roots 
$\bt\not\in \{\al,\frac{1}{2}\al,2\al\}$ with $(\al,\bt)(\al,\al)>0$, $\bt+n\al$ is a root for all positive integers $n$ unless $\al$ and $\bt$ are both positive or negative, $\bt$ is of finite type and norm zero and $\al-\bt$ is a root.)
\item
Let $\al$ and $\bt$ be orthogonal positive (resp. negative) roots of
infinite type or of norm 0.  If $x\in G_{\al}$ and $[x, G_{-\gamma}]=0$ for
all roots $\gamma$ for which $0< \left| \gamma (h) \right| <
\left| \al(h) \right|$, then $[x,G_{\bt}]=0$.
\end{enumerate}
Then $G$ is a direct sum of a generalized Kac-Moody superalgebra, finite
dimensional simple Lie superalgebras 
and affine Lie superalgebras. 
\end{thth1}

\section{The vertex algebra of the compactified superstring}

In this section we construct a vertex algebra that represents the Fock space of a chiral superstring moving on a 10 dimensional torus and describe its symmetries. 
 
\subsection{The construction of the vertex algebra}

Now we construct the vertex algebra of the compactified superstring using the lattice construction described in section \ref{lc}. 

The 10 dimensional superstring has 10 bosonic and 10 fermionic degrees of freedom. We represent the bosonic degrees by a 10 dimensional even lattice and the fermionic degrees in bosonized form by the weightlattice of the Lie algebra $D_5\,$. The conformal ghosts and superghosts give 3 additional lattice dimensions. One of them carries a negative metric and must be described together with the fermionic lattice (cf. \cite{LT}). 

Hence we define
\be L = L^X \oplus L^{\psi,\phi} \oplus L^{\chi,\sigma} \ee
where $L^X$ is a 10 dimensional even Lorentzian lattice, 
$L^{\psi,\phi}=\{ (x_1,\ldots,x_6)\in {\mathbb R}^{5,1}\, | \, \mbox{all }
x_i\in {\mathbb Z}\, \mbox{ or all } x_i\in {\mathbb Z}+\frac{1}{2}\}$
and $L^{\chi,\sigma}={\mathbb Z}^2$
 
The lattice $L^{\psi,\phi}$ is a 6 dimensional rational lattice. 
The elements $(x_1,\ldots,x_6)$ with integral $x_i$ and $\sum x_i$ even form an even sublattice of $L^{\psi,\phi}$. The quotient of $L^{\psi,\phi}$ by this lattice is isomorphic to ${\mathbb Z}_2 \times {\mathbb Z}_2$ and will be written $\Gamma(L^{\psi,\phi})=\{ 0,V,S,C\}$. We have the coset decomposition $L^{\psi,\phi}=L^{\psi,\phi}_0 \cup L^{\psi,\phi}_V \cup L^{\psi,\phi}_S \cup L^{\psi,\phi}_C$ with
\beann 
L^{\psi,\phi}_V &=& (0,\ldots,0,1) + L^{\psi,\phi}_0 \\
L^{\psi,\phi}_S &=& ({\textstyle \frac{1}{2}},\ldots,
                     {\textstyle \frac{1}{2}},{\textstyle \frac{1}{2}})
                                                  + L^{\psi,\phi}_0 \\
L^{\psi,\phi}_C &=& ({\textstyle \frac{1}{2}},\ldots,
                     {\textstyle \frac{1}{2}},-{\textstyle \frac{1}{2}})
                                                  + L^{\psi,\phi}_0 
\eeann
We can write the conjugacy classes of $L^{\psi,\phi}$ also in terms of the classes of $D_5^*$ and $D_1^*$. The lattice $D_n^*$ has the decomposition 
$0(D_n^*)\cup V(D_n^*)\cup  S(D_n^*)\cup C(D_n^*)$ where 
$0(D_n^*)$ is the even lattice with elements 
$(x_1,\ldots,x_n)\in {\mathbb Z}^n$ 
such that $\sum x_i$ is even, and
\beann
V(D_n^*) &=& (0,\ldots,0,1)+0(D_n^*) \\ 
S(D_n^*) &=& ({\textstyle\frac{1}{2}},\ldots,{\textstyle\frac{1}{2}},{\textstyle\frac{1}{2}})+0(D_n^*) \\
C(D_n^*) &=& ({\textstyle\frac{1}{2}},\ldots,{\textstyle\frac{1}{2}},-{\textstyle\frac{1}{2}})+0(D_n^*) \, . 
\eeann
The quotient of $D_n^*$ by $0(D_n^*)$ is ${\mathbb Z}_2\times {\mathbb Z}_2$ for $n$ even and ${\mathbb Z}_4$ for $n$ odd. 
With this notation we can write 
\[ \renewcommand{\arraystretch}{1.3} \begin{array}{c|c}
L^{\psi,\phi} & (D_5^*,D_1^*) \\ \hline 
       0      & (0,0)\cup (V,V) \\
       V      & (0,V)\cup (V,0) \\
       S      & (S,S)\cup (C,C) \\
       C      & (C,S)\cup (C,S) 
\end{array} \]
An element of $L^{\psi,\phi}$ has even normsquare if it is in $L^{\psi,\phi}_0$ and odd normsquare else. 
The vectors $\phi^i,\, i=1,\ldots,5$ with $1$ in the i-th position and zero else and 
$s^6=({\textstyle \frac{1}{2}},\ldots,{\textstyle \frac{1}{2}})$ form a 
$\mathbb Z$-basis of $L^{\psi,\phi}$
The sublattice $E_{5,1}=L^{\psi,\phi}_0 \cup L^{\psi,\phi}_S$ is an integral unimodular Lorentzian lattice.

The lattice $L^{\chi,\sigma}={\mathbb Z}^2$ can be decomposed as 
$L^{\chi,\sigma}=L^{\chi,\sigma}_0 \cup L^{\chi,\sigma}_1$ into the elements of even and odd norm. The quotient $\Gamma(L^{\chi,\sigma})=\{ 0, 1\}$ is isomorphic to ${\mathbb Z}_2$. A $\mathbb Z$-basis of $L^{\chi,\sigma}$ is given by the elements $\chi=(1,0)$ and $\sigma=(0,1)$. 

Now the sum $L^X \oplus L^{\psi,\phi}_0 \oplus L^{\chi,\sigma}_0$ is an even sublattice of $L$. The quotient $\Gamma(L)$ of $L$ by this sublattice is $\Gamma(L^{\psi,\phi})\times \Gamma(L^{\chi,\sigma})$. 

To construct a vertex algebra from $L$ we need a bilinear map $\eta$ from $\Gamma(L) \times \Gamma(L)$ to ${\mathbb C}^*$. There is a number of such maps but we need a special one to assure that the modes of the physical bosons and fermions that we will define have the correct commutation and anticommutation properties.

The following definition of $\eta$ turns out to be appropriate.

\[ \renewcommand{\arraystretch}{1.3} \begin{array}{c|cccccccc}
      & (0,0) & (V,0) & (S,0) & (C,0) & (0,1) & (V,1) & (S,1) & (C,1) \\ \hline
(0,0) &   1   &   1   &   1   &   1   &   1   &   1   &   1   &   1   \\
(V,0) &   1   &  -1   &   y   &  -y   &   -1  &   1   &   -y  &   y   \\
(S,0) &   1   &  -y   &  -1   &   y   &   -1  &   y   &   1   &   -y  \\
(C,0) &   1   &   y   &  -y   &  -1   &   1   &   y   &   -y  &   -1  \\
(0,1) &   1   &  -1   &  -1   &   1   &   -1  &   1   &   1   &   -1  \\
(V,1) &   1   &   1   &  -y   &  -y   &   1   &   1   &   -y  &   -y  \\
(S,1) &   1   &   y   &   1   &   y   &   1   &   y   &   1   &    y  \\  
(C,1) &   1   &   -y  &   y   &   -1  &   -1  &   y   &   -y  &    1
\end{array} \]

$y$ can be chosen as $+1$ or $-1$. The column denotes the first and the row the second argument of $\eta$. Note that $\eta$ is not symmetric, e.g. 
$\eta ((S,0),(V,0))=y$ and $\eta ((V,0),(S,0))=-y$. One can check that this map is bilinear. It clearly satisfies condition (\ref{nc}).

Let $\{\al^1,\dots,\al^{10}\}$ be a basis of $L^X$. Then 
$\{ \al^1,\dots,\al^{10},\phi^1,\dots,\phi^5,s^6,\chi,\sigma \}$ is an 
ordered basis of $L$. As 2-cocycle we choose the one constructed from this basis as described above with the $a$'s being as follows. Let $i\leq j$. Then define $a_{ii}=(-1)^{\frac{1}{2}(\al^i,\al^i)}$ for $i=1,\ldots,10$ and $a_{ij}=1$
in the other cases. This definition will give nice formulas.

We define the vertex algebra $V$ of the compactified superstring as the vertex algebra constructed from the lattice $L$ with $\eta$ and $\varepsilon$ given as above. Hence
\[
V =          V_{(0,0)} \oplus V_{(V,0)} \oplus V_{(S,0)} \oplus V_{(C,0)} 
      \oplus V_{(0,1)} \oplus V_{(V,1)} \oplus V_{(S,1)} \oplus V_{(C,1)}\, .
\]
$V$ can be decomposed into the Neveu-Schwarz and Ramond sector given by 
\[ V^{NS}=V_{(0,0)} \oplus V_{(V,0)} \oplus V_{(0,1)} \oplus V_{(V,1)} \]
and   
\[ V^R=V_{(S,0)} \oplus V_{(C,0)} \oplus V_{(S,1)} \oplus V_{(C,1)} \, . \]

The vector space 
\be V^{GSO} =  
V_{(0,0)} \oplus  V_{(S,1)} \oplus V_{(0,1)} \oplus V_{(S,0)} \ee
describes the Fock space of a GSO-projected superstring. $V^{GSO}$ is a vertex superalgebra with even part $V_{(0,0)} \oplus  V_{(S,1)}$ and odd part  $V_{(0,1)} \oplus V_{(S,0)}$.

\subsection{The symmetries of the compactified superstring}

In this section we show that $V$ carries representations of the $N=1$ and $N=2$ superconformal algebras and use them to construct the BRST operator.  

First we will consider the matter sector.

Let $\{ x^1,\dots,x^{10} \}$ be a basis of ${\mathbb R}\otimes_{\mathbb Z} L^X$ with $(x^{\mu},x^{\nu})= g^{\mu \nu}$ where $g^{\mu \nu}$ is the diagonal matrix with entries $g^{1 1}=\ldots=g^{9 9}=1$ and $g^{10\, 10}=-1$. The inverse of $g^{\mu \nu}$ is denoted by $g_{\mu \nu}$. These matrices can be used to raise and lower indices. The bosons $x^{\mu}(-1)$ satisfy 
$x^{\mu}(-1)_0x^{\nu}(-1)=0,\, x^{\mu}(-1)_1x^{\nu}(-1)=g^{\mu \nu}$ and $x^{\mu}(-1)_nx^{\nu}(-1)=0$ for $n\geq 2$. By (\ref{dpe}) this is equivalent to the operator product expansion
\be x^{\mu}(-1)(z)\,x^{\nu}(-1)= g^{\mu \nu} z^{-2} + \ldots \ee
where the dots denote nonsingular terms.

Next we define 10 complex fermions 
\[ \Psi^{\pm i} = e^{\pm \phi^i} \quad i=1,\ldots,5 \, .  \]
Since the complex fermions are elements of $V_{(V,0)}$ their vertex operators have integral expansions when acting on the Neveu-Schwarz sector and half-integral expansions when acting on the Ramond sector. They satisfy
\be \Psi^{i}(z)\Psi^{-j}= \delta^{ij} z^{-1}+\ldots \, . \ee
The 10 real fermions 
\[ \renewcommand{\arraystretch}{1.3} \begin{array}{lcl c lcl}
 \psi^1 & = &  \frac{1}{\sqrt 2} (\Psi^1 + \Psi^{-1}) & \qquad &
 \psi^2 & = &  \frac{i}{\sqrt 2} (\Psi^1 - \Psi^{-1}) \\
        & \vdots &   & & & \vdots &                  \\ 
 \psi^7 & = &  \frac{1}{\sqrt 2} (\Psi^4 + \Psi^{-4}) & &
 \psi^8 & = &  \frac{i}{\sqrt 2} (\Psi^4 - \Psi^{-4}) \\

 \psi^9 & = &  \frac{1}{\sqrt 2} (\Psi^5 + \Psi^{-5}) & & 
 \psi^{10} & = & \frac{1}{\sqrt 2} (\Psi^5 - \Psi^{-5})
\end{array} \]
have the operator product expansion  
\be  \psi^{\mu}(z)\psi^{\nu} = g^{\mu \nu} z^{-1} + \ldots \, . \ee
Using the commutator formula this implies 
\be  \{ \psi^{\mu}_m,\psi^{\nu}_n \} = g^{\mu \nu} \delta_{m+n+1} \ee
where $m$ and $n$ are either both integral or both half-integral depending on wether the $\psi$'s act on the Neveu-Schwarz or Ramond sector.

Define the supercurrent and the energy momentum tensor of the matter sector by 
\be \tau^M = g_{\mu \nu} x^{\mu}(-1)_{-1}\psi^{\nu} 
           = x_{\nu}(-1)_{-1}\psi^{\nu}              \ee
and
\be
  \omega^M = {\textstyle \frac{1}{2}} \tau^M_{\; 0}\,\tau^M 
           = {\textstyle \frac{1}{2}} x_{\nu}(-1)_{-1}x^{\nu}(-1)
             +{\textstyle \frac{1}{2}}  (D\psi_{\nu})_{-1}\psi^{\nu} \, . 
\ee
Note that we assume the summation convention in these expressions. 
We remark that the energy momentum tensor can also be written as
\[ \omega^M={\textstyle \frac{1}{2}} x_{\nu}(-1)_{-1}x^{\nu}(-1) 
+ {\textstyle \frac{1}{2}}\phi^1(-1)_{-1}\phi^1(-1) + \ldots 
+ {\textstyle \frac{1}{2}}\phi^5(-1)_{-1}\phi^5(-1) \, .\]
It is easy to calculate the following operator product expansions
\bea
\omega^M(z)x^{\mu}(-1) &=& x^{\mu}(-1)z^{-2} + Dx^{\mu}(-1)z^{-1} + \dots \\
\omega^M(z)\psi^{\mu}  &=& {\textstyle \frac{1}{2}}\psi^{\mu}z^{-2} 
                           + D\psi^{\mu}z^{-1} + \dots 
\eea
and
\bea
\tau^M(z)x^{\mu}(-1) &=& \psi^{\mu}z^{-2} + D\psi^{\mu}z^{-1} + \dots \\
\tau^M(z)\psi^{\mu}  &=& x^{\mu}(-1)z^{-1} + \dots \, .
\eea

To obtain a representation of the $N=2$ superconformal algebra we define 
\be    j^M = \phi^1(-1) + \ldots + \phi^5(-1) \, . \ee
and decompose the supercurrent 
\be   \tau^M = \tau^{M+} + \tau^{M-} \ee
with
\[ \tau^{M+} = \sum_{i=1}^5 h^i(-1)_{-1}\Psi^i \, ,\qquad
   \tau^{M-} = \sum_{i=1}^5 h^{-i}(-1)_{-1}\Psi^{-i} \]
where $h^1 =  (x^1+ix^2)/{\sqrt 2},\ldots,
       h^4 =  (x^7+ix^8)/{\sqrt 2},
       h^5 =  (x^9-x^{10})/{\sqrt 2}$ 
and $ h^{-1} = (x^1-ix^2)/{\sqrt 2},\ldots,
      h^{-4} = (x^7-ix^8)/{\sqrt 2},
      h^{-5} = (x^9+x^{10})/{\sqrt 2}$. 
These elements satisfy $(h^i,h^j)=(h^{-i},h^{-j})=0$ and $(h^{i},h^{-j})=\delta^{ij}$.

Using the associativity formula one proves the following 
\begin{p1}
The above fields satisfy the operator product expansions
\[ \renewcommand{\arraystretch}{1.3} 
\begin{array}{rlccclclcl}
\omega^M(z) & \hspace{-0.3cm}\omega^M &=& {\textstyle \frac{15}{2}} z^{-4} 
                   &+& 2 \omega^M z^{-2} &+& D\omega^M z^{-1} &+& \ldots  \\
\omega^M(z) & \hspace{-0.3cm}\tau^{M\pm} &=& 
                              {\textstyle \frac{3}{2}} \tau^{M\pm} z^{-2}
                               &+& D\tau^{M\pm} z^{-1} &+& \dots & & \\
\omega^M(z) & \hspace{-0.3cm}j^M  &=& j^M z^{-2} 
                               &+& Dj^M z^{-1} &+& \dots & & \\
j^M(z) & \hspace{-0.3cm}\tau^{M\pm} &=& \pm \tau^{M\pm} z^{-1} &+& \dots & & 
                                                                       & &  \\
j^M(z) & \hspace{-0.3cm}j^M &=& 5 z^{-2} &+& \dots & &   & &  
\end{array}\]
\[ \renewcommand{\arraystretch}{1.3} 
\begin{array}{rlclclclcl}
\tau^{M+}(z) & \hspace{-0.3cm}\tau^{M-} &=& 5 z^{-3} &+& j^M z^{-2}  
                           &+& (\omega^M+ {\textstyle \frac{1}{2}}Dj^M) z^{-1}
                           &+& \dots \\
\tau^{M-}(z) & \hspace{-0.3cm}\tau^{M+} &=& 5 z^{-3} &-& j^M z^{-2}  
                           &+& (\omega^M - {\textstyle \frac{1}{2}}Dj^M) z^{-1}
                           &+& \dots 
\end{array}\]
The expansions $\tau^{M+}(z)\tau^{M+}$ and $\tau^{M-}(z)\tau^{M-}$ contain only nonsingular terms.
\end{p1}
Now define operators $L^M_{m}=\omega^M_{m+1}$ and 
$G^{M\pm}_{m}=\tau^{M\pm}_{m+\frac{1}{2}}$. Then the commutator formula and the above operator product expansions imply 
\[ \renewcommand{\arraystretch}{1.3} 
\begin{array}{lcl}
 {[} L^M_{m}, L^M_{n} {]} & = & (m-n)L^M_{m+n} + \frac{1}{12}m(m^2-1)
                                                  \delta_{m+n}c^M \\
 {[} L^M_{m}, G^{M\pm}_{n} {]} & = & (\frac{1}{2}m-n) G^{M\pm}_{m+n} \\
 {[} L^M_{m}, j^M_{n} {]}         & = & -n j^M_{m+n} \\
 {[} j^M_m, G^{M\pm}_{n} {]}   & = & \pm G^{M\pm}_{m+n} \\
 {[} j^M_m,j^M_n {]}    & = & \frac{1}{3}m \delta_{m+n}c^M \\
 \{ G^{M+}_m,G^{M-}_n\} & = & L^M_{m+n} + \frac{1}{2}(m-n) j^M_{m+n} 
                        + \frac{1}{6}(m^2-\frac{1}{4}) \delta_{m+n}c^M \\
 \{ G^{M+}_m,G^{M+}_n\} & = & \{ G^{M-}_m,G^{M-}_n\} = 0
\end{array} \]
where the central charge is $c^M=15$. The expansions of the fields $\tau^{M\pm}(z)$ and $j^M(z)$ depand on where they act. However the corresponding algebras are all isomorphic to the $N=2$ superconformal algebra. The isomorphism is called spectral flow. 

The expansions for the $\tau^{M\pm}$ give 
\[ \tau^M(z)\tau^M = 10 z^{-3} + 2\omega^{M} z^{-1} + \ldots \] 
so that the operators $G^{M}_{m}=\tau^M_{m+\frac{1}{2}}=G^{M+}_m+G^{M-}_m$ satisfy 
\[ \renewcommand{\arraystretch}{1.3}
\begin{array}{lcl}
 {[} L^M_{m}, L^M_{n} {]} & = & (m-n)L^M_{m+n} + \frac{1}{12}m(m^2-1)
                                                  \delta_{m+n}c^M \\
 {[} L^M_{m}, G^M_{n} {]} & = & (\frac{1}{2}m-n)G^M_{n} \\
  \{ G^{M}_m,G^{M}_n \}   & = & 2 L^M_{m+n} 
                       + \frac{1}{3} (m^2-\frac{1}{4}) \delta_{m+n}c^M
\end{array} \]
where $m$ and $n$ are half-integral in the Neveu-Schwarz sector and integral in the Ramond sector. The corresponding 2 superalgebras are called Neveu-Schwarz and Ramond algebra.

Now we show that there is another representation of the $N=2$ superconformal algebra on $V$ coming from the ghost fields.

Let $\phi=(0,\ldots,0,1)\in L^{\psi,\phi}$. 
Note that $\varepsilon (\phi, \phi) =-1$ and 
$\varepsilon (\phi-\chi,-\phi+\chi)=\varepsilon (\sigma, \sigma)=1$.
Define the ghost fields
\be b=e^{-\sigma}\, , \quad  c=e^{\sigma} \ee
and
\be \beta = \chi(-1)_{-1}e^{-\phi+\chi} \, , \quad \gamma = e^{\phi-\chi} \, . \ee
Then
\[ c(z)b = z^{-1}+\ldots \quad \mbox{and} \quad 
   \gamma(z)\beta = z^{-1}+\ldots \: .              \]
The other expansions between these fields contain only non singular terms.
This implies that 
\[ \{ c_m ,b_n \} = \delta_{m+n+1} \quad \mbox{and} \quad
    [ \gamma_m,\beta_n ] = \delta_{m+n+1} \]
are the only nontrivial commutation relations. 

Now we define the supercurrent and the energy momentum tensor of the ghost sector by
\be \tau^{Gh}=c_{-1}(D\beta)+{\textstyle\frac{3}{2}}(Dc)_{-1}\beta
              -2b_{-1}\gamma \ee
and 
\be \omega^{Gh} = {\textstyle \frac{1}{2}} \tau^{Gh}_{\; 0}\,\tau^{Gh} 
             = 2(Dc)_{-1}b - (Db)_{-1}c  
               - {\textstyle \frac{3}{2}}(D\gamma)_{-1}\beta
               - {\textstyle \frac{1}{2}}(D\beta)_{-1}\gamma \, .\ee
The energy momentum tensor can be written as 
\be 
\omega^{Gh} = \omega^{\phi} + \omega^{\chi} + \omega^{\sigma} 
\ee
where the elements
\beann
\omega^{\phi}   &=& -{\textstyle \frac{1}{2}}\phi(-1)_{-1}\phi(-1)
                         -\phi(-2) \\  
\omega^{\chi}   &=& {\textstyle \frac{1}{2}}\chi(-1)_{-1}\chi(-1) 
                         +{\textstyle \frac{1}{2}}\chi(-2)          \\
\omega^{\sigma} &=& {\textstyle \frac{1}{2}}\sigma(-1)_{-1}\sigma(-1)
                         + {\textstyle \frac{3}{2}}\sigma(-2)
\eeann 
generate representations of the Virasoro algebra on $V$ of central charges 
$13,-2$ and $-26$.

To get a representation of the $N=2$ superconformal algebra we define
\be j^{Gh} = -3\phi(-1) + 2 \sigma(-1)\, . \ee
and decompose the supercurrent as
\be \tau^{Gh}=\tau^{Gh+}+\tau^{Gh-} \ee
with 
\be
    \tau^{Gh+}=c_{-1}(D\beta)+{\textstyle\frac{3}{2}}(Dc)_{-1}\beta \, , \quad
    \tau^{Gh-}=-2b_{-1}\gamma \, .
\ee
As for the matter fields we get by using the associativity formula 
\begin{pp1} 
The above fields satisfy the following operator product expansions
\[ \renewcommand{\arraystretch}{1.3}
\begin{array}{rlccclclcl}
\omega^{Gh}(z) & \hspace{-0.3cm}\omega^{Gh} 
                        &=& -{\textstyle \frac{15}{2}} z^{-4} 
                        &+& 2 \omega^{Gh} z^{-2} 
                        &+& D\omega^{Gh} z^{-1} &+& \ldots  \\
\omega^{Gh}(z) & \hspace{-0.3cm}\tau^{Gh\pm} 
                        &=& {\textstyle \frac{3}{2}} \tau^{Gh\pm} z^{-2}
                        &+& D\tau^{Gh\pm} z^{-1} &+& \dots & & \\
\omega^{Gh}(z) & \hspace{-0.3cm}j^{Gh}  
                        &=& j^M z^{-2} 
                        &+& Dj^{Gh} z^{-1} &+& \dots & & \\
j^{Gh}(z) & \hspace{-0.3cm}\tau^{Gh\pm} 
                        &=& \pm \tau^{Gh\pm} z^{-1} &+& \dots & &   & &  \\
j^{Gh}(z) & \hspace{-0.3cm}j^{Gh} 
                        &=& -5 z^{-2} &+& \dots & &   & &  
\end{array}\]
\[ \renewcommand{\arraystretch}{1.3} 
\begin{array}{rlclclclcl}
\tau^{Gh+}(z) & \hspace{-0.3cm}\tau^{Gh-} 
                     &=& -5 z^{-3} &+& j^{Gh} z^{-2}  
                     &+& (\omega^{Gh}+ {\textstyle \frac{1}{2}}Dj^{Gh}) z^{-1}
                     &+& \dots \\
\tau^{{Gh}-}(z) & \hspace{-0.3cm}\tau^{Gh+} 
                     &=& -5 z^{-3} &-& j^{Gh} z^{-2}  
                     &+& (\omega^{Gh} - {\textstyle \frac{1}{2}}Dj^{Gh}) z^{-1}
                     &+& \dots 
\end{array}\]
The expansions $\tau^{Gh+}(z)\tau^{Gh+}$ and $\tau^{Gh-}(z)\tau^{Gh-}$ contain only nonsingular terms.
\end{pp1}
Hence the operators $L^{Gh}_{n}=\omega^{Gh}_{n+1}, 
G^{Gh\pm}_{n}=\tau^{Gh\pm}_{n+\frac{1}{2}}$ and $j^{Gh}_n$ give a representation of the $N=2$ superconformal algebra of central charge $c^{Gh}=-15$. 

As before the opertors $L^{Gh}_{n}$ and $G^{Gh}_{n}=G^{Gh+}_{n}+G^{Gh-}_{n}$
give a representation of the Neveu-Schwarz resp. Ramond algebra in the Neveu-Schwarz resp. Ramond sector of central charge $c^{Gh}=-15$.

We give some expansions of the energy momentum tensor and the supercurrent of the ghost sector. They can be used to prove the next proposition.
\[ \renewcommand{\arraystretch}{1.3} 
\begin{array}{lcccccl}
\omega^{Gh}(z)b &=& 2b z^{-2} &+& Db z^{-1} &+& \ldots \\
\omega^{Gh}(z)c &=& -c z^{-2} &+& Dc z^{-1} &+& \ldots \\
\omega^{Gh}(z)\beta &=& \frac{3}{2}\beta z^{-2} &+& D\beta z^{-1} &+&\ldots \\
\omega^{Gh}(z)\gamma &=& -\frac{1}{2}\gamma z^{-2} &+& D\gamma z^{-1} 
                                       &+&\ldots \\
                     & &   & &   & & \\ 
\tau^{Gh}(z)b &=& -\frac{3}{2}\beta z^{-2} &-& \frac{1}{2}D\beta z^{-1} 
                                           &+& \ldots \\
\tau^{Gh}(z)c      &=& &-& 2\gamma z^{-1} &+& \ldots \\
\tau^{Gh}(z)\beta  &=& &-& 2b z^{-1}       &+&\ldots \\
\tau^{Gh}(z)\gamma &=& c z^{-2} &-& \frac{1}{2}Dcz^{-1} &+&\ldots 
\end{array}\]

We define $\omega=\omega^{M}+\omega^{Gh}$ and $\tau=\tau^{M}+\tau^{Gh}$.
Then $\omega$ is a Virasoro element, in particular $D=L_{-1}$. The operators $L_{n}=\omega_{n+1}$ and $G_{n}=\tau_{n+\frac{1}{2}}$ give a representation of the Neveu-Schwarz and the Ramond algebra of central charge zero.

Let 
\be j^{BRST}= c_{-1}(\omega^M+{\textstyle \frac{1}{2}}\omega^{Gh})
              +\gamma_{-1}(\tau^M+{\textstyle \frac{1}{2}}\tau^{Gh})\, . \ee
The zero mode of $j^{BRST}$ is called BRST operator and denoted by $Q$.
Using the associativity formula and the above expansions we get

\pagebreak[3]
\begin{pp1}
The BRST operator transforms the above fields as
\bea 
Qb      &=& \omega \\
Qc      &=& -c_{-2}c-\gamma_{-1}\gamma \\
Q\beta  &=& \tau \\
Q\gamma &=& \gamma_{-2}c - {\textstyle \frac{1}{2}}c_{-2}\gamma \\
Q\omega &=& 0 \\
Q\tau   &=& 0 \, .
\eea
\end{pp1}
These formulas imply 
\bea 
\{Q,b_n\} &=& \omega_n \\
\{Q,c_n\} &=& -(c_{-2}c+\gamma_{-1}\gamma)_n \label{qc} \\
{[}Q,\beta_n{]}  &=& \tau_n \\
{[}Q,\gamma_n{]} &=& (\gamma_{-2}c - {\textstyle \frac{1}{2}}c_{-2}\gamma)_n 
                                     \label{qg}\\
{[}Q,\omega_n{]}   &=& 0 \\
\{ Q,\tau_n \} &=& 0 \, .
\eea

Equations (\ref{qc}) and (\ref{qg}) together with the explicit expressions for $Q\omega^{Gh}$ and $Q\tau^{Gh}$ show that $Qj^{BRST} = Dv$ with
$  v= \frac{1}{4}c_{-1}\gamma_{-1}\tau^M 
      + \frac{1}{4}\gamma_{-1}\gamma_{-1}\gamma_{-1}\beta$ 
      $- \frac{1}{2}\gamma_{-1}\gamma_{-1}c_{-1}b
      + \frac{3}{2}c_{-3}c + \frac{3}{2}\gamma_{-2}\gamma$ .
Hence \be Q^2=0 \, . \ee

We derive a decomposition of $Q$ that will be useful for calculations in the next section. Write $\omega^{Gh}=\omega^{b,c} + \omega^{\beta,\gamma}$  where 
\[ \renewcommand{\arraystretch}{1.3}
\begin{array}{c}
\omega^{b,c}=2(Dc)_{-1}b - (Db)_{-1}c=\omega^{\sigma} \\
\omega^{\beta,\gamma}
 =-\frac{3}{2}(D\gamma)_{-1}\beta - \frac{1}{2}(D\beta)_{-1}\gamma 
 =\omega^{\phi}+\omega^{\chi} \, .
\end{array} \]
Then
$j^{BRST}=c_{-1}(\omega^M+\omega^{\beta,\gamma}
                                      +\frac{1}{2}\omega^{b,c})
          + \gamma_{-1}\tau^M - \gamma_{-1}\gamma_{-1}b 
          + \frac{3}{4}D(\beta_{-1}\gamma_{-1}c)$ 
implies that 
\be Q = Q_0 + Q_1 +Q_2 \ee
with
\beann
Q_0 &=& \left(c_{-1}(\omega^M+\omega^{\beta,\gamma}
                     +{\textstyle \frac{1}{2}}\omega^{b,c})\right)_0 \\
Q_1 &=& \left(\gamma_{-1}\tau^M\right)_0 \\
Q_2 &=& -\left(\gamma_{-1}\gamma_{-1}b\right)_0 \, .
\eeann 
The following lemma is easy to prove.
\begin{lp1} \label{ecl}
Let $v$ be a state satisfying
\beann
 (L^M_0+L^{\beta,\gamma}_0)v &=& v \\
 (L^M_n+L^{\beta,\gamma}_n)v &=& 0 \quad\mbox{for}\quad n\geq 1 \\
  L^{b,c}_n v &=& 0  \quad\mbox{for}\quad n\geq -1 \, .
\eeann
Then \[ Q_0v =D(c_{-1}v) \, .\]
\end{lp1}

Finally we construct operators that induce gradings on $V$. Let $j=a\phi(-1)+b\chi(-1)+c\sigma(-1)$. The zero mode 
$j_0=a\phi(0)+b\chi(0)+c\sigma(0)$ of $j$ measures the ghost charge of a state.
The BRST element $j^{BRST}$ is an eigenvector of $j_0$ iff $a+b+c=0$.
We define 
\[ \renewcommand{\arraystretch}{1.3}
\begin{array}{lcl}
 j^{N} &=& -\chi(-1)+\sigma(-1) \\
 j^{P} &=& -\phi(-1)  +\chi(-1) \, .
\end{array} \]
The corresponding zero modes are called ghost number and ghost picture 
operator. 
A state has integral ghost picture if it is in the Neveu-Schwarz sector and half-integral ghost picture if it is in the Ramond sector. 
The elements of $V_{(.., 0)}$ have even and the elements of $V_{(.., 1)}$ odd ghost number.

We list the $L_0$-eigenvalue and the ghost number and ghost picture of some 
states.
\[ \renewcommand{\arraystretch}{1.3}
\begin{array}{l|ccc}
            & L_0                & j^{N}_0 & j^{P}_0 \\ \hline
e^{n\phi}   & -\frac{1}{2}n(n+2) & 0       & n       \\
e^{n\chi}   & \frac{1}{2}n(n-1)  & -n      & n       \\
e^{n\sigma} & \frac{1}{2}n(n-3)  & n       & 0       \\
j^{BRST}    &  1                 & 1       & 0
\end{array} \]
The last row implies 
\[ [j^{N}_0,Q]=Q \quad\mbox{ and }\quad [j^{P}_0,Q]=0 \, . \]
Hence the BRST operator increases the ghost number by one and leaves the ghost picture unchanged.

\section{The Lie algebra of physical states}

In this section we construct the Lie superalgebra of physical states of the 10 dimensional superstring moving on a torus and study its properties.
Section \ref{mrs} contains the main results of this paper. 

\subsection{The small algebra}

Here we define a subalgebra of the vertex algebra of the compactified superstring which is called small algebra. This is necessary to define the picture changing operator which is essential for the construction of the Lie superalgebra of physical states. We also construct a bilinear form on the small algebra.

Define the elements $\xi=e^{\chi}$ and $\eta=e^{-\chi}$. We have 
\[ \renewcommand{\arraystretch}{1.3}
\begin{array}{lcrll}
\xi(z)\eta   &=& z^{-1}  & +& \ldots \\
D \xi(z)\eta &=& -z^{-2} & +& \ldots \, .
\end{array} \]
The expansions $D\xi(z)D\xi$ and $\eta(z)\eta$ have no singular terms so that
\[ \renewcommand{\arraystretch}{1.3}
\begin{array}{lcl}
\{ (D\xi)_m,\eta_n \}   &=& -m\delta_{m+n} \\
\{ (D\xi)_m,(D\xi)_n \} &=& 0 \\
\{ \eta_m, \eta_n \}    &=& 0 \, .
\end{array}\]

We define the small algebra $V_S$ as the subalgebra of $V$ generated 
by the states $e^{\gamma},D\xi,\eta,b,c$ with 
$\gamma\in L^X \oplus L^{\psi,\phi}$. Note that $\beta$ and $\gamma$ are in the small algebra.

Let
\[ X= Q\xi = (Q_0 + Q_1 + Q_2)\xi 
           = c_{-1}(D\xi) + \tau^M_{-1}e^{\phi} + D(\eta_{-1}e^{2\phi}_{-1}b)
                             + (D\eta)_{-1}e^{2\phi}_{-1}b \, .
\]
This is an element in $V_{(0,0)}$. The mode $X_{-1}=(Q\xi)_{-1}=\{Q,\xi_{-1}\}$ is an endomorphism of the small algebra called picture changing operator.

\begin{p1}
The picture changing operator has the following properties
\[ \renewcommand{\arraystretch}{1.3}
\begin{array}{lcl}
{[} X_{-1}, Q  \, ] &=& 0 \\
{[} X_{-1}, L_0 \,] &=& 0 \\
{[} X_{-1}, b_n \,] &=& -(D\xi)_{n-1} \, .
\end{array}  \]
In particular $[ X_{-1}, b_1 ] =0$.
\end{p1}
{\em Proof:} The first relation is clear. The second one follows from $L_0\xi=0$. $X_k=Q\xi_k+\xi_kQ$ and $\xi_kb=0$ for $k\geq 0$ imply that for $k\geq 0$
\beann
X_kb &=& \xi_kQb \\
     &=& \xi_k\omega \\
     &=& -[ L_{-2},\xi_k]1 \\
     &=& -\sum_{m\geq 0}{-1 \choose m} (L_{m-1}\xi)_{-1+k-m} 1 \\
     &=& -(D\xi)_{k-1}1 
\eeann
since $L_{m-1}\xi=0$ for $m\geq 1$. By the commutator formula we have
\[ [X_{-1},b_n\,]=\sum_{k\geq 0}{-1 \choose k} (X_kb)_{-1+n-k} =(X_0b)_{n-1}
               =-(D\xi)_{n-1} \, . \]
\hspace*{\fill} $\Box$\\

$X$ has ghost number $0$ and ghost picture $1$. This implies
\[ [j^{N}_0,X_{-1}]=0 \quad\mbox{ and }\quad [j^{P}_0,X_{-1}]=X_{-1} \, . \]

It is easy to see that the eigenvalues of $L_0$ on $V^{GSO}$ are all integral. Moreover $L_1$ acts locally nilpotent on $V^{GSO}$. The same holds for the subalgebra $V^{GSO}_S=V_S\cap V^{GSO}$ so that we can apply the results of section \ref{ibf} to construct a bilinear form on this subalgebra.
Since $e^{3\sigma-2\phi}$ is not in the image of $L_1$ we can define an invariant bilinear form $(\, ,\, )$ on $V^{GSO}_S$ by setting 
$(e^{3\sigma-2\phi},1)=1$. It is easy to prove
\begin{pp1}
$(\, ,\, )$ is nondegenerate on $V^{GSO}_S$. 
$(e^{\gamma+n\sigma},e^{\gamma'+n'\sigma})$ with 
$\gamma,\gamma'\in L^X \oplus L^{\psi,\phi}$ is nonzero if $\gamma+\gamma'=-2\phi$ and $n+n'=3$
and zero else.
\end{pp1}

We calculate the adjoints of some operators. $b$ and $c$ have $L_0$-eigenvalue $2$ and $-1$ and are annihilated by $L_1$. Hence 
\be b_n^* = b_{2-n} \quad\mbox{ and }\quad c_n^*=-c_{-4-n}\, . \ee
From $L_1j^{BRST}=\frac{3}{2}Dc,\, 
L_1^2j^{BRST}=-3c$ and $L_1^nj^{BRST}$ $=0$ for $n\geq 2$ we get
\be Q^*=-Q \,. \ee
Finally $L_0X=L_1X=0$ gives 
\be X_{-1}^* = X_{-1} \, . \ee 

\subsection{The cohomology spaces}

We now construct the cohomology spaces that represent the physical states of the compactified superstring and determine their dimensions.

The space 
\[ B = V^{GSO}_S \cap \mbox{Ker } b_1 \]
is graded by the lattice $L^X$. We denote the homogeneous subspaces by 
$B(\al)$. The operators $L_0,\, j_0^N$ and $j_0^P$ are simultaneously diagonizable on $B(\al)$ so that we can write 
\[ B(\al)=\bigoplus B(\al)_{p,n}^k \]
where $k,n$ and $p$ denote the eigenvalues of $L_0,\, j_0^N$ and $j_0^P$.
We define spaces $C(\al)_{p,n}=B(\al)_{p,n}^0,\, C(\al)=\bigoplus C(\al)_{p,n}$ and $C=\bigoplus C(\al)$.

The following proposition can be proved by analyzing the $L_0$-contributions of the ghost sector. 
\begin{pp1}
Let $\al$ be in $L^X$. Then
\begin{enumerate}
\item[]
The space $B(\al)_{p,n}^k$ is finite dimensional.
\item[]
The space $B(\al)_{-1}^k=\bigoplus_{n\in \mathbb Z} B(\al)_{-1,n}^k$ is
finite dimensional.
\item[]
The space $B(\al)_{-\frac{1}{2}}^k
=\bigoplus_{n\in \mathbb Z} B(\al)_{-\frac{1}{2},n}^k$ is in general not 
finite dimensional.
For example the state
\[ e^{\lambda+(n-\frac{1}{2})\phi}_{-1}\eta_{-n}\ldots\eta_{-1}c\, , \]
where $\lambda=(\pm\frac{1}{2},\ldots,\pm\frac{1}{2})$ is in $C(D_5^*)$ if $n$ is even and in $S(D_5^*)$ if $n$ is 
odd, is in the kernel of $b_1$, has $L_0$-eigenvalue zero, ghost picture $-\frac{1}{2}$ and ghost number $n$.
\item[]
The spaces $C(\al)_{-1,1}$ and $C(\al)_{-\frac{1}{2},1}$ are zero for 
$(\al,\al)>0$.
\end{enumerate}
Now let $\al$ be of norm zero. Then
\begin{enumerate}
\item[]
$C(\al)_{-1,1}$ is 10 dimensional and spanned by the states 
$\psi^{\mu}_{-1}e^{-\phi}_{-1}c_{-1}e^{\al}$. 
\item[]
$C(\al)_{-1,0}$ is spanned by $e^{-2\phi}_{-1}(D\xi)_{-1}c_{-1}e^{\al}$.
\item[]
The 16 states 
$e^{\lambda_{\dot{\beta}}-\frac{1}{2}\phi}_{-1}c_{-1}e^{\al}$ 
with 
$\lambda_{\dot{\beta}}=(\pm\frac{1}{2},\ldots,\pm\frac{1}{2})\in C(D_5^*)$ 
form a basis of $C(\al)_{-\frac{1}{2},1}$.
\item[]
The space $C(\al)_{-\frac{1}{2},0}$ is zero.
\item[]
The 16 states 
$e^{\lambda_{\beta}-\frac{3}{2}\phi}_{-1}c_{-1}e^{\al}$ 
with $\lambda_{\beta}=(\pm\frac{1}{2},\ldots,\pm\frac{1}{2})\in S(D_5^*)$ 
are a basis of $C(\al)_{-\frac{3}{2},1}$.
\item[]
$C(\al)_{-\frac{3}{2},0}$ is spanned by the 16 states
$e^{\lambda_{\dot{\beta}}-\frac{5}{2}\phi}_{-1}(D\xi)_{-1}c_{-1}e^{\al}$
where 
$\lambda_{\dot{\beta}}=(\pm\frac{1}{2},\ldots,\pm\frac{1}{2})$ is in 
$C(D_5^*)$.
\end{enumerate}
\end{pp1}

Since $b_1^*=b_1$ and $\{b_1,c_{-2}\}=1$ the bilinear form on $C$ induced from $(\, ,\, )$ is zero. We can construct a nontrivial bilinear form on 
$C$ by defining 
\[ (u,v)_C=(c_{-2}u,v)\, . \]
Then
\begin{pp1} \label{bf1}
$(\, ,\, )_C$ pairs nondegenerately $C(\al)_{p,m}$ and 
$C(\bt)_{q,n}$ where $\al+\bt=0,\, p+q=-2$ and $m+n=2$. In particular
$(\, ,\, )_C$ is nondegenerate on $C$.

Furthermore $(Qu,v)_C$ is up to a sign equal to $(u,Qv)_C$.
\end{pp1}
{\em Proof:} $\{b_1,c_{-2}\}=1$ implies that $c_{-2}$ is an injective on the kernel of $b_1$. If $u\in C(\al)_{p,m}$ is nonzero then there is an element $v$ in $V^{GSO}_S$ such that $(c_{-2}u,v)$ is nonzero. Since $[L_0,c_{-2}]=0$ this element has $L_0$-eigenvalue zero. Clearly it has ghost number $2-m$, ghost picture $-2-q$ and momentum $-\al$. The commutator $\{b_1,c_{-2}\}=1$ shows that $v$ is in the kernel of $b_1$. 

For $u,v\in C$ 
\beann
(Qu,v)_C &=& (c_{-2}Qu,v) \\
                &=& (c_{-2}Qb_1c_{-2}u,v) \\
                &=& -(c_{-2}b_1Qc_{-2}u,v) \\
                &=& -(Qc_{-2}u,v) \\
                &\sim & (c_{-2}u,Qv) \\
                &=& (u,Qv)_C
\eeann
\hspace*{\fill} $\Box$\\
Since $Q^2=0$ the sequence 
\[ \ldots\stackrel{Q}{\longrightarrow} C(\al)_{p,n-1} 
         \stackrel{Q}{\longrightarrow} C(\al)_{p,n}
         \stackrel{Q}{\longrightarrow} C(\al)_{p,n+1}
         \stackrel{Q}{\longrightarrow} \ldots \]
defines a complex. We denote the cohomology groups by $H(\al)_{p,n}$. Define  
the spaces $H(\al)=\bigoplus H(\al)_{p,n}$ and $H=\bigoplus H(\al)$. 
We investigate the structure of the space $H$.

Proposition \ref{bf1} and the fact that the spaces $C(\al)_{p,n}$ are 
finite dimensional give 
\begin{pp1}
$(\, ,\, )_C$ induces a nondegenerate bilinear form 
$(\, ,\, )_H$ on $H$.
\end{pp1}

As in \cite{LZ1} one can prove
\begin{thp1}
Let $\al\neq 0$. Then 
\[ H(\al)_{-1,n}=0 \quad \mbox{for } n\neq 1\,. \]
\end{thp1}

We use this result and the Euler-Poincar\'{e} principle to calculate the dimension of $H(\al)_{-1,1}$ for $\al\neq 0$.
Define numbers $c(n)$ by
\bea
 \sum c(n)q^n &=& 8 \prod_{m=1}^{\infty}\left(\frac{1+q^m}{1-q^m}\right)^8 \\
              &=& 8+128q+1152q^2+7680q^3+42112q^4+\ldots \, .
\eea
The asymptotic behaviour of the $c(n)$ is given by 
\be c(n)\sim {\textstyle \frac{1}{2}}
                         n^{-\frac{11}{4}}e^{2\pi\sqrt{2n}} \, . \ee
For $n\geq 10$ the right hand side gives a good approximation for $c(n)$.

We have
\begin{thp1}
Let $\al\neq 0$. Then 
\[ \mbox{dim}\, H(\al)_{-1,1}=c({\textstyle -\frac{1}{2}}\al^2)\,. \]
\end{thp1}
{\em Proof:} Let $\al\neq 0$. 
Since the space $C_{-1}(\al)=\bigoplus_{n\in \mathbb Z} C_{-1,n}(\al)$ is finite dimensional we can apply the Euler-Poincar\'{e} principle to get
\[ - \mbox{dim}\, H(\al)_{-1,1}=  \sum(-1)^n \mbox{dim}\,C(\al)_{-1,n}\, . \]
The right hand side is the constant term in
\[ \mbox{tr}\, (-1)^{j_0^N}q^{L_0} \]
where the trace has to be evaluated on 
$B(\al)_{-1}=\bigoplus_{k,n \in \mathbb Z}B(\al)_{-1,n}^k$.
Since the spaces $B(\al)_{-1}^k$ are finite dimensional $\mbox{tr}\, (-1)^{j_0^N}q^{L_0}$ is well defined as a formal Laurent series in $q$.
By looking at a basis of $B(\al)_{-1}$ we find 
\begin{enumerate}
\item[] 
The 10 bosonic coordinates $x^{\mu}$ give a contribution
\[    q^{\frac{1}{2}\al^2} \varphi(q)^{-10} \]
to the trace where $\varphi(q) = \prod_{n=1}^{\infty} (1-q^n)$.
\item[] 
The $(\psi,\phi)$-$(D\xi,\eta)$ system gives a contribution 
\[ \varphi(q)^{-7}\sum (-1)^m q^{\frac{1}{2}\lambda^2}
      q^{-\frac{1}{2}p^2+\frac{1}{2}m(m+1)+\frac{1}{2}} \]
where the sum extends over $(\lambda,-p-1)\in L_0^{(\psi,\phi)}$ and 
$m\in \mathbb Z$ with $m\geq |p|$. 
Note that the sum is a well defined power series in $q$.
\item[]
The $(b,c)$ system contributes a factor
\[ -q^{-1}\varphi(q)^2 \]
to the trace.
\end{enumerate}
Hence the trace is given by
\[ -q^{\frac{1}{2}\al^2-1}\varphi(q)^{-15}
      \sum (-1)^m q^{\frac{1}{2}\lambda^2}
           q^{-\frac{1}{2}p^2+\frac{1}{2}m(m+1)+\frac{1}{2}}\, \]
A technical calculation 
shows that this is equal to
\[ - 8 q^{\frac{1}{2}\al^2} \frac{\varphi(q^2)^{8}}{\varphi(q)^{16}}\,. \]
The constant term in this expression is $-c(-\frac{1}{2}\al)$.
This proves the theorem. \hspace*{\fill} $\Box$\\

We will see that for $\al\neq 0$ the spaces $H(\al)_{-1,1}$ and 
$H(\al)_{-\frac{1}{2},1}$ are related by a supersymmetry transformation so that they have the same dimension. 

The picture changing operator commutes with $Q$ and therefore defines a map on $H$. The following proposition (cf. \cite{BZ}) describes an important property of this map.
\begin{pp1} \label{piso}
For $\al\neq 0$ the map $X_{-1}: H(\al)_{p,n}\rightarrow H(\al)_{p+1,n}$ is an isomorphism.
\end{pp1}
{\em Proof:} $\al=\al_{\nu}x^{\nu}$ has a nonzero component $\al_{\mu}$. 
The element $\tilde{P}^{\mu}=\psi^{\mu}_{-1}e^{-\phi}$ has ghost number $0$ and ghost picture $-1$. The zero mode $\tilde{P}^{\mu}_0$ is the momentum operator in the $-1$ picture. It satisfies 
$[Q,\tilde{P}^{\mu}_0]=(Q\tilde{P}^{\mu})_0=(Q_0\tilde{P}^{\mu})_0
=(D(c_{-1}\tilde{P}^{\mu}))_0=0$. 
Hence $\tilde{P}^{\mu}_0$ defines a map from $H(\al)_{p,n}$ to $H(\al)_{p-1,n}$. One can show that $\tilde{P}^{\mu}_0$ is a scalar multiple of the inverse of $X_{-1}$ \linebreak[3] (cf. \cite{BZ}). \hspace*{\fill} $\Box$\\
 
Hence for nonzero momentum the spaces $H(\al)_{p,1}$ with integral or half-integral ghost picture $p$ are all isomorphic. They represent the bosonic resp. fermionic physical states of momentum $\al$. The ghost pictures $-1$ and 
$-\frac{1}{2}$ are called canonical ghost pictures because they give the easiest description of the physical states. The fermionic physical states also admit a simple description in the $-\frac{3}{2}$ picture since the contribution of the $\phi$-charge to the $L_0$-value in this picture is the same as in the $-\frac{1}{2}$ picture. 

Before we can determine the cohomology spaces in the canonical ghost pictures for $\al^2=0$ explicitly we have to introduce some more notations.

Let 
\[ S^{\dot{\al}}=e^{\lambda_{\dot{\al}}-\frac{1}{2}\phi} \quad\mbox{and}\quad 
   S^{\al}=e^{\lambda_{\al}-\frac{3}{2}\phi}                              \]
where 
$\lambda_{\dot{\al}}=(\pm\frac{1}{2},\ldots,\pm\frac{1}{2})$ is a conjugate spinor of $D_5^*$, i.e. has an odd number of $-$ signs, and
$\lambda_{\al}=(\pm\frac{1}{2},\ldots,\pm\frac{1}{2})$ is a spinor of $D_5^*$ and has an even number of $-$ signs.

We define $\Gamma$-matrices $\Gamma^{\mu}$ by
\bea
(\psi^{\mu}_{-1}e^{-\phi})(z)S^{\dot{\al}} &=& \frac{1}{\sqrt{2}}
{\Gamma^{\mu\, \dot{\al}}}_{\bt}S^{\bt}z^{-1} + \ldots \\
(\psi^{\mu}_{-1}e^{\phi})(z)S^{\al} &=& \frac{1}{\sqrt{2}}
{\Gamma^{\mu\, \al}}_{\dot{\bt}}S^{\dot{\bt}}z^{1} + \ldots \, .
\eea
The upper index denotes the row and the lower index the column.
\begin{pp1}
We have
\bea
(\psi^{\mu}_{-1}e^{\phi})_{-2}(\psi^{\nu}_{-1}e^{-\phi})_{0}S^{\dot{\al}} &=& 
\psi^{\mu}_{-\frac{1}{2}}\psi^{\nu}_{-\frac{1}{2}}S^{\dot{\al}} \\
(\psi^{\mu}_{-1}e^{-\phi})_{0}(\psi^{\nu}_{-1}e^{\phi})_{-2}S^{\al} &=& 
\psi^{\mu}_{-\frac{1}{2}}\psi^{\nu}_{-\frac{1}{2}}S^{\al}\, .
\eea
\end{pp1}
{\em Proof:}
\beann
(\psi^{\mu}_{-1}e^{\phi})_{-2}(\psi^{\nu}_{-1}e^{-\phi})_0S^{\dot{\al}} &=&
(\psi^{\mu}_{-1}e^{\phi})_{-2} 
\psi^{\nu}_{-\frac{1}{2}}e^{-\phi}_{-\frac{1}{2}}S^{\dot{\al}} \\
 &=& (\psi^{\mu}_{-\frac{3}{2}}e^{\phi}_{-\frac{3}{2}}
       -e^{\phi}_{-\frac{5}{2}}\psi^{\mu}_{-\frac{1}{2}})
      \psi^{\nu}_{-\frac{1}{2}}e^{-\phi}_{-\frac{1}{2}}S^{\dot{\al}} \\
 &=& -e^{\phi}_{-\frac{5}{2}}\psi^{\mu}_{-\frac{1}{2}}
        \psi^{\nu}_{-\frac{1}{2}}e^{-\phi}_{-\frac{1}{2}}S^{\dot{\al}}\\
 &=& -\psi^{\mu}_{-\frac{1}{2}}\psi^{\nu}_{-\frac{1}{2}}
     e^{\phi}_{-\frac{5}{2}}e^{-\phi}_{-\frac{1}{2}}S^{\dot{\al}}\\
 &=& \psi^{\mu}_{-\frac{1}{2}}\psi^{\nu}_{-\frac{1}{2}} S^{\dot{\al}} 
\eeann
We have used (\ref{bisc}) repeatedly.
Note that
$e^{\phi}_{-\frac{5}{2}}e^{-\phi}_{-\frac{1}{2}}S^{\dot{\al}}=-S^{\dot{\al}}$
since 
\[ \varepsilon(-\phi,\lambda_{\dot{\al}}-{\textstyle \frac{1}{2}}\phi)
       \varepsilon(\phi,\lambda_{\dot{\al}}-{\textstyle \frac{3}{2}}\phi)
       =\varepsilon(\phi,\lambda_{\dot{\al}}-{\textstyle \frac{1}{2}}\phi)^{-1}
        \varepsilon(\phi,\lambda_{\dot{\al}}-{\textstyle \frac{1}{2}}\phi)
        \varepsilon(\phi,\phi)^{-1}=-1\,. \]
The proof of the other equation is analogous.  \hspace*{\fill} $\Box$\\

The commutator $ \{ \psi^{\mu}_{-\frac{1}{2}},\psi^{\nu}_{-\frac{1}{2}} \}
                 = g^{\mu \nu}$ implies 
\be \{ \Gamma^{\mu},\Gamma^{\nu}\}=2g^{\mu \nu} \ee
so that the $\Gamma$-matrices give a Weyl representation of the Clifford algebra $C(9,1)$. The matrix $\Gamma^{11}=\Gamma^{1}\ldots \Gamma^{10}$ anticommutes with the $\Gamma^{\mu}$ and satisfies $(\Gamma^{11})^2=1$. More precisely $\Gamma^{11}$ is diagonal and with diagonal entries $\pm(1,\ldots,1,-1,\ldots,-1)$.
The operators $\frac{1}{2}(1\pm \Gamma^{11})$ are projection operators on the spinors and conjugate spinors.

We define a charge conjugation matrix $C$ by
\bea
S^{\al}(z)S^{\dot{\bt}} &=& C^{\al \dot{\bt}}e^{-2\phi}z^{-2}+\ldots \\
S^{\dot{\al}}(z)S^{\bt} &=& C^{\dot{\al} \bt}e^{-2\phi}z^{-2}+\ldots  \, .
\eea
The definition implies that $C$ is antisymmetric and invertible. We have
\begin{pp1} \label{susyv}
The spinor fields have the operator product expansions
\beann
S^{\dot{\al}}(z)S^{\dot{\bt}} &=& \frac{1}{\sqrt{2}}
                                  (\Gamma_{\mu}C)^{\dot{\al} \dot{\bt}}
                                  \psi^{\mu}_{-1}e^{-\phi} z^{-1} + \ldots \\
S^{\al}(z)S^{\bt} &=&            \frac{1}{\sqrt{2}}      
                                 (\Gamma_{\mu}C)^{\al \bt}
                                  \psi^{\mu}_{-1}e^{-3\phi} z^{-3} + \ldots 
\eeann
\end{pp1}
{\em Proof:} 
It is easy to see that the products $S^{\dot{\al}}_nS^{\dot{\bt}}$ vanish for $n\geq 1$. We make the ansatz $S^{\dot{\al}}_{\, 0}S^{\dot{\bt}}=
M_{\mu}^{\dot{\al} \dot{\bt}}\psi^{\mu}_{-1}e^{-\phi}$. Equation (\ref{bisc}) with 
$n=0,\, k=1$ and $m=0$ gives
\[ ((\psi^{\mu}_{-1}e^{-\phi})_0S^{\dot{\al}})_1S^{\dot{\bt}}
= (\psi^{\mu}_{-1}e^{-\phi})_1 S^{\dot{\al}}_0S^{\dot{\bt}} \, .\]
Inserting the definitions and using again equation (\ref{bisc}) we get
\[ \frac{1}{\sqrt{2}}
     {\Gamma^{\mu\, \dot{\al}}}_{\gamma} C^{\gamma \dot{\bt}}e^{-2\phi} 
   = M_{\nu}^{\dot{\al} \dot{\bt}}
     (\psi^{\mu}_{-1}e^{-\phi})_1(\psi^{\nu}_{-1}e^{-\phi})  
   = M_{\nu}^{\dot{\al} \dot{\bt}} g^{\mu \nu} e^{-2\phi}\, . \]
This implies the first equation.\\
The second equation is proved in the same way.
\hspace*{\fill} $\Box$\\
Clearly $\Gamma^{\mu} C$ is symmetric and 
$C^{-1}\Gamma^{\mu}C=-{\Gamma^{\mu}}^T$.

Now we describe the cohomology spaces for $\al^2=0$.
\pagebreak[3]
\begin{pp1}
Let $\al\in L^X$ with $\al^2=0$ and $\al\neq 0$. Then
\begin{enumerate}
\item[]
$H(\al)_{-1,1}$ is generated by the massless vectors
$|\xi,\al\ra = \xi_{\mu}\psi^{\mu}_{-1}e^{-\phi}_{-1}c_{-1}e^{\al}$
satisfying $(\xi,\al)=\xi_{\mu}\al^{\mu}=0$. The state 
$\al_{\mu}\psi^{\mu}_{-1}e^{-\phi}_{-1}c_{-1}e^{\al}$ is trivial in the 
cohomology. $H(\al)_{-1,1}$ has dimension $8$.
\item[]
$H(\al)_{-\frac{1}{2},1}$ also has dimension $8$ and is spanned by the spinors $|u,-\frac{1}{2},\al \ra = u_{\dot{\bt}}S^{\dot{\bt}}_{-1}c_{-1}e^{\al}$ satisfying the massless Dirac equation 
$\al_{\mu}{\Gamma^{\mu\, \dot{\bt}}}_{\gamma} u_{\dot{\bt}}=0$.
\item[]
$H(\al)_{-\frac{3}{2},1}$ is generated by the massless spinors $|u,-\frac{3}{2},\al \ra = u_{\bt}S^{\bt}_{-1}c_{-1}e^{\al}$. The states with
$\al_{\mu}{\Gamma^{\mu \, \bt}}_{\dot{\delta}}u_{\bt}=0$ are BRST-trivial.
$H(\al)_{-\frac{3}{2},1}$ also has dimension $8$.
\end{enumerate}
Let $\al=0$. Then
\begin{enumerate}
\item[]
$H(\al)_{-1,1}$ has dimension $10$. The vectors 
$P^{\mu}=-\psi^{\mu}_{-1}e^{-\phi}_{-1}c$ form a basis of this space.
\item[]
$H(\al)_{-\frac{1}{2},1}$ has dimension $16$ and a basis 
$Q^{\dot{\al}}=S^{\dot{\al}}_{-1}c$.
\item[]
$H(\al)_{-\frac{3}{2},1}$ has dimension $16$ and a basis 
$Q^{\al}=S^{\al}_{-1}c$.
\end{enumerate}
\end{pp1} 
{\em Proof:} We can write $\xi_{\mu}\psi^{\mu}_{-1}e^{-\phi}_{-1}c_{-1}e^{\al}=c_{-1}w$ with $w=\xi_{\mu}\psi^{\mu}_{-1}e^{-\phi}_{-1}e^{\al}$. 
Then by Lemma \ref{ecl}
\[ Qc_{-1}w = \{Q,c_{-1}\}w-c_{-1}Qw = -c_{-1}Q_1w \, .\]
Using equation (\ref{bisc}) with $n=k=-1$ and $m=1$ we get 
\[ Q_1w =(\gamma_{-1}\tau^M)_0 w 
        = \sum_{j\geq 0}(\gamma_{-2-j}\tau^M_{j+1}+\tau^M_{-j}\gamma_{-1+j})w 
        = \gamma_{-2}\tau^M_1 w  \, .\]
Now $\tau^M_1w = (\xi,\al) e^{-\phi}_{-1}e^{\al}$. Calculating 
$Qe^{-2\phi}_{-1}(D\xi)_{-1}c_{-1}e^{\al}$ with Lemma \ref{ecl} shows that  
$\al_{\mu}\psi^{\mu}_{-1}e^{-\phi}_{-1}c_{-1}e^{\al}$ is trivial in cohomology.
This proves the first statement. \\
The massless spinor 
$u_{\dot{\bt}}S^{\dot{\bt}}_{-1}c_{-1}e^{\al}$ is in the kernel of $Q$ iff 
$Q_1w=0$ where $w=u_{\dot{\bt}}S^{\dot{\bt}}_{-1}e^{\al}$. Equation 
(\ref{bisc}) gives 
\[ Q_1 w = \sum_{j\geq 0}( \gamma_{-\frac{3}{2}-j}\tau^M_{j+\frac{1}{2}}
                        +\tau^M_{-\frac{1}{2}-j}\gamma_{-\frac{1}{2}+j})w 
         = \gamma_{-\frac{3}{2}}\tau^M_{\frac{1}{2}} w  \]
and 
\beann 
\tau^M_{\frac{1}{2}} w 
&=& (x_{\nu}(-1)_{-1}\psi^{\nu})_{\frac{1}{2}} w \\
&=& x_{\nu}(-1)_0\psi^{\nu}_{-\frac{1}{2}} w \\
&=& \al_{\nu}u_{\dot{\bt}}e^{\al}_{-1}\psi^{\nu}_{-\frac{1}{2}} S^{\dot{\bt}}\\
&=& -\al_{\nu}u_{\dot{\bt}}e^{\al}_{-1}\psi^{\nu}_{-\frac{1}{2}} 
    e^{\phi}_{-\frac{5}{2}}e^{-\phi}_{-\frac{1}{2}}S^{\dot{\bt}} \\
&=& \al_{\nu}u_{\dot{\bt}}e^{\al}_{-1}e^{\phi}_{-\frac{5}{2}}
    \psi^{\nu}_{-\frac{1}{2}}e^{-\phi}_{-\frac{1}{2}}S^{\dot{\bt}} \\
&=& \frac{1}{\sqrt{2}}\al_{\nu}u_{\dot{\bt}}{\Gamma^{\nu\, \dot{\bt}}}_{\gamma}
    e^{\phi}_{-\frac{5}{2}}S^{\gamma}_{-1}e^{\al}
\eeann
where we have used 
$S^{\dot{\bt}}=-e^{\phi}_{-\frac{5}{2}}e^{-\phi}_{-\frac{1}{2}}S^{\dot{\bt}}$. 
Hence $Q_1w=0$ iff 
$\al_{\nu}{\Gamma^{\nu\,\dot{\bt}}}_{\gamma}u_{\dot{\bt}}=0$.
Since $C(\al)_{-\frac{1}{2},0}=0$ this proves the second statement.\\
It is easy to see that the elements in $C(\al)_{-\frac{3}{2},1}$ are in the kernel of $Q$. Thus we have to calculate the image of $C(\al)_{-\frac{3}{2},0}$ 
under $Q$. We have 
\[ Qu_{\dot{\delta}}e^{\lambda_{\dot{\delta}}-\frac{5}{2}\phi}_{-1}(D\xi)_{-1}
   c_{-1}e^{\al}
   = Qc_{-1}w = -c_{-1}Q_1w \]
with $w=u_{\dot{\delta}}e^{\lambda_{\dot{\delta}}-\frac{5}{2}\phi}_{-1}
      (D\xi)_{-1}e^{\al}$.
As above 
\[ Q_1w = (\gamma_{-1}\tau^M)_0w  
        = \sum_{j\geq 0}( \gamma_{-\frac{3}{2}-j}\tau^M_{j+\frac{1}{2}})w
        = \gamma_{-\frac{3}{2}}\tau^M_{\frac{1}{2}} w \, .\] 
Now
\beann
\tau^M_{\frac{1}{2}} w
&=& \al_{\mu}u_{\dot{\delta}}(D\xi)_{-1}e^{\al}_{-1}
    \psi^{\mu}_{-\frac{1}{2}}e^{\lambda_{\dot{\delta}}-\frac{5}{2}\phi} \\
&=& \al_{\mu}u_{\dot{\delta}}(D\xi)_{-1}e^{\al}_{-1}
    \psi^{\mu}_{-\frac{1}{2}}
    (e^{-2\phi}_0e^{\lambda_{\dot{\delta}}-\frac{1}{2}\phi}) \\
&=& \al_{\mu}u_{\dot{\delta}}(D\xi)_{-1}e^{\al}_{-1}e^{-2\phi}_0
    \psi^{\mu}_{-\frac{1}{2}}e^{\lambda_{\dot{\delta}}-\frac{1}{2}\phi} \\
&=& \frac{1}{\sqrt{2}}\al_{\mu}u_{\dot{\delta}}
    {\Gamma^{\mu \, \dot{\delta}}}_{\gamma}
    (D\xi)_{-1}e^{\al}_{-1}e^{-2\phi}_0
    e^{\phi}_{-\frac{5}{2}}
    e^{\lambda_{\gamma}-\frac{3}{2}\phi},
\eeann
where we have used 
$\varepsilon(-2\phi,\lambda_{\dot{\delta}}-{\textstyle\frac{1}{2}}\phi) =1$
and 
\beann
\lefteqn{Q_1w} \\ 
&=& \frac{1}{\sqrt{2}}\al_{\mu}u_{\dot{\delta}}
    {\Gamma^{\mu \, \dot{\delta}}}_{\gamma}e^{\al}_{-1}
    \gamma_{-\frac{3}{2}}(D\xi)_{-1}e^{-2\phi}_0
    e^{\phi}_{-\frac{5}{2}}e^{\lambda_{\gamma}-\frac{3}{2}\phi}  \\
&=& \frac{1}{\sqrt{2}}\al_{\mu}u_{\dot{\delta}}
    {\Gamma^{\mu \, \dot{\delta}}}_{\gamma}e^{\al}_{-1}
    (e^{\phi}_{-1}\eta)_{-\frac{3}{2}}(D\xi)_{-1}e^{-2\phi}_0
    e^{\phi}_{-\frac{5}{2}}e^{\lambda_{\gamma}-\frac{3}{2}\phi}  \\
&=& \frac{1}{\sqrt{2}}\al_{\mu}u_{\dot{\delta}}
    {\Gamma^{\mu \, \dot{\delta}}}_{\gamma}e^{\al}_{-1}
    \sum_{j\geq 0}(e^{\phi}_{-\frac{9}{2}-j} \eta_{2+j}
                   -\eta_{1-j}e^{\phi}_{-\frac{7}{2}+j}) 
    (D\xi)_{-1}e^{-2\phi}_0
    e^{\phi}_{-\frac{5}{2}}e^{\lambda_{\gamma}-\frac{3}{2}\phi} \\
&=& \frac{1}{\sqrt{2}}\al_{\mu}u_{\dot{\delta}}
    {\Gamma^{\mu \, \dot{\delta}}}_{\gamma}e^{\al}_{-1}
    \eta_1 (D\xi)_{-1} e^{\phi}_{-\frac{7}{2}}e^{-2\phi}_0
    e^{\phi}_{-\frac{5}{2}}e^{\lambda_{\gamma}-\frac{3}{2}\phi} \\
&=& -\frac{1}{\sqrt{2}}\al_{\mu}u_{\dot{\delta}}
    {\Gamma^{\mu \, \dot{\delta}}}_{\gamma}e^{\al}_{-1}
    e^{\lambda_{\gamma}-\frac{3}{2}\phi} 
\eeann
since
\pagebreak[3] 
\beann
\lefteqn{
   \varepsilon(\phi,\lambda_{\gamma}-{\textstyle \frac{3}{2}}\phi)
   \varepsilon(-2\phi,\lambda_{\gamma}-{\textstyle \frac{1}{2}}\phi)
   \varepsilon(\phi,\lambda_{\gamma}-{\textstyle \frac{5}{2}}\phi)} \\
&=& \varepsilon(\phi,\lambda_{\gamma}-{\textstyle \frac{3}{2}}\phi)
    \varepsilon(-2\phi,\lambda_{\gamma}-{\textstyle \frac{1}{2}}\phi) 
    \varepsilon(\phi,\lambda_{\gamma}-{\textstyle \frac{3}{2}}\phi)
    \varepsilon(\phi,-\phi) \\
&=& \varepsilon(2\phi,\lambda_{\gamma}-{\textstyle \frac{3}{2}}\phi)
    \varepsilon(-2\phi,\lambda_{\gamma}-{\textstyle \frac{1}{2}}\phi) 
    \varepsilon(\phi,\phi)^{-1} \\
&=& \varepsilon(2\phi,-\phi) \varepsilon(\phi,\phi)^{-1} \\
&=& -1 \, . 
\eeann
Hence 
\[ Q u_{\dot{\delta}}e^{\lambda_{\dot{\delta}}-\frac{5}{2}\phi}_{-1}
   (D\xi)_{-1}c_{-1}e^{\al}
   = -\frac{1}{\sqrt{2}}\al_{\mu}{\Gamma^{\mu \, \dot{\delta}}}_{\gamma}
     u_{\dot{\delta}} S^{\gamma}_{-1}c_{-1}e^{\al}\, . \]
Thus $v_{\bt}S^{\bt}_{-1}c_{-1}e^{\al}$ is BRST-trivial iff 
$v_{\bt}=\al_{\mu}{\Gamma^{\mu \, \dot{\delta}}}_{\bt}u_{\dot{\delta}}$
which is equivalent to 
$\al_{\mu}{\Gamma^{\mu \, \bt}}_{\dot{\delta}}v_{\bt}=0$.
This proves the third statement. \\
The rest of the proposition is now clear.
\hspace*{\fill} $\Box$\\

We conclude this section with some remarks on the picture changing operator. 
First we give two examples. 
The massless vector is given in the $0$ picture by
\be  X_{-1}|\xi,\al \ra = 
\left( \xi_{\mu}x^{\mu}(-1)_{-1}+\al_{\nu}\xi_{\mu}
  \psi^{\nu}_{-1}\psi^{\mu}_{-1}\right) e^{\al}_{-1}c 
  + \xi_{\mu}\psi^{\mu}_{-1}e^{\al}_{-1}\gamma \, . \ee
For $|u,-\frac{3}{2},\al \ra = u_{\bt}S^{\bt}_{-1}e^{\al}_{-1}c$ we have
\be X_{-1} |u,{\textstyle -\frac{3}{2}},\al \ra 
   = \frac{1}{\sqrt 2} \al_{\mu}{\Gamma^{\mu\,\bt}}_{\dot{\gamma}} 
     u_{\bt}S^{\dot{\gamma}}_{-1}c_{-1}e^{\al} \, . \ee
Although $\mbox{dim }H(\al)_{-\frac{1}{2},1}=H(\al)_{-\frac{3}{2},1}=16$ for 
$\al=0$ the last equation implies
\begin{pp1} \label{ncgp}
Let $\al=0$. Then $X_{-1}H(\al)_{-\frac{3}{2},1}=0$.
\end{pp1}
This shows that Proposition \ref{piso} does not hold for zero momentum.

\subsection{The Lie algebra and the Lie bracket}

We define a product on $V^{GSO}_S$ that induces a product on the cohomology space $H$ (cf. \cite{LZ2}) and construct a Lie bracket on the physical states using the picture changing operator.

For $u,v$ in $V^{GSO}_S$ define
\be \{u,v\}=(-1)^{|u|}(b_0u)_0v \ee
where $|u|$ denotes the parity of $u$ with respect to $V^{GSO}_S$. Then
\begin{pp1}\label{qder}
For $u,v$ in $V^{GSO}_S$ we have
\beann
(-1)^{|u|}\{ u,v \} &=& b_1(u_{-1}v)-(b_1u)_{-1}v-(-1)^{|u|} u_{-1}(b_1v) \\
Q\{ u,v \} &=& \{Qu,v\}+(-1)^{|u|+1}\{u,Qv\} \, .
\eeann
\end{pp1}
Since $b_1^2=0$ the first equation implies that the bracket closes on $B$.
It is easy to see that it also closes on $C$. The second equation shows that it projects down to the cohomology $H$. Hence $\{\, ,\, \}$ defines a product
$H(\al)_{p,m}\times H(\bt)_{q,n}\rightarrow H(\al+\beta)_{p+q,m+n-1}$.
The product satisfies (cf. \cite{LZ2})
\begin{pp1} \label{bsym}
Let $u,v$ and $w$ be representatives of elements in $H$. Then
\[ \renewcommand{\arraystretch}{1.5}
\begin{array}{c}
 \{ u,v \}+ (-1)^{(|u|+1)(|v|+1)}\{ v,u \} = 0 \\
  (-1)^{(|u|+1)(|w|+1)}\{u,\{v,w\}\} +  (-1)^{(|v|+1)(|u|+1)}\{v,\{w,u\}\} \\
                  + (-1)^{(|w|+1)(|v|+1)}\{w,\{u,v\}\} =0 
\end{array} \]
\end{pp1}
Now we show that the bracket invariant under picture changing.
\begin{pp1} \label{pinv}
Let $u,v\in H$. Then
\[ X_{-1}\{u,v\}=\{X_{-1}u,v\}=\{u,X_{-1}v\} \]
\end{pp1}
{\em Proof:}
From Proposition (\ref{qder}) we get
\[ (-1)^{|u|}X_{-1} \{ u,v \} - (-1)^{|u|}\{ u, X_{-1}v\} 
              = b_1 [X_{-1},u_{-1}]v \, . \]
Now
\[ [X_{-1},u_{-1}]=\sum_{k\geq 0} {-1 \choose k} (X_ku)_{-k-2}
                  =\sum_{k\geq 0} {-1 \choose k} (Q\xi_ku)_{-k-2} \]
and 
\[ (Q\xi_ku)_{-k-2}= \left( j_0^{BRST} (\xi_ku) \right)_{-k-2}
                   = Q(\xi_ku)_{-k-2}\pm (\xi_ku)_{-k-2}Q \]
imply
\[
[X_{-1},u_{-1}]v = \sum_{k\geq 0} {-1 \choose k}Q(\xi_ku)_{-k-2}v \\
                 = \sum_{k\geq 0} \frac{(-1)^{k+1}}{k+1} Q 
                     \left( (D\xi)_{k+1}u \right)_{-k-2}v \, .
\]
Note that $\left((D\xi)_{k+1}u\right)_{-k-2}v$ is an element of $V^{GSO}_S$
with $L_0$-eigenvalue zero. $b_1Q=L_0-Qb_1$ shows then that 
$b_1 [X_{-1},u_{-1}]v$ is in the image of $Q$. This proves 
$X_{-1}\{u,v\}=\{u,X_{-1}v\}$. The Proposition now follows from the symmetry of the product. \hspace*{\fill} $\Box$\\

Lian and Zuckerman have defined in \cite{LZ2} another product $u\cdot v=u_{-1}v$. By the derivation property of zero modes this product defines a product on the cohomology. It satisfies $u\cdot v =(-1)^{|u| |v|}v\cdot u$. Then as above $X_{-1}(u\cdot v) -u\cdot (X_{-1}v)=[X_{-1},u_{-1}]v$ is zero in cohomology so that this product is also invariant under $X_{-1}$, i.e. 
$X_{-1}(u\cdot v)=u\cdot (X_{-1}v)=(X_{-1}u)\cdot v$.

Define ${\tilde G}(\al)={\tilde G}(\al)_0\oplus {\tilde G}(\al)_1$ where 
${\tilde G}(\al)_0 = H(\al)_{-1,1}$ and ${\tilde G}(\al)_1 = H(\al)_{-\frac{1}{2},1}$.
Let \[ {\tilde G} = \bigoplus_{\al\in L^X} {\tilde G}(\al) \]
and define a ${\mathbb Z}_2$-grading on ${\tilde G}$ by setting ${\tilde G}_0=\bigoplus {\tilde G}(\al)_0$
and ${\tilde G}_1=\bigoplus {\tilde G}(\al)_1$.

We define a product on ${\tilde G}$ by
\[ [u,v] = X_{-1} \{u,v\} \quad \mbox{ if $u$ or $v$ is in ${\tilde G}_0$} \]
and 
\[ [u,v] = \{u,v\}  \quad \mbox{ if $u$ and $v$ are in ${\tilde G}_1$}\, .\]
Note that the parity of $u$ as an element of ${\tilde G}$ is equal to $|u|+1$ 
where $|u|$ denotes the parity of $u$ as an element in $V^{GSO}_S$. 
By Propositions \ref{bsym} and \ref{pinv} we have
\begin{pp1}
With this product ${\tilde G}$ is a Lie superalgebra.
\end{pp1}
Define $G$ as ${\tilde G}$ without ${\tilde G}(0)_1$, i.e. without the supersymmetry charges $Q^{\dot{\al}}$. Then Proposition \ref{ncgp} shows that $G$ is a subalgebra of ${\tilde G}$. We will see that $G$ is in contrast to ${\tilde G}$ a generalized Kac-Moody superalgebra. Since $G$ represents the physical states of the compactified superstring in the canonical ghost pictures we call $G$ the Lie superalgebra of physical states. The even part of $G$ contains the bosonic and the odd part the fermionic physical states. We have already seen that $G$ has no real roots. This corresponds to the fact that the superstring has no tachyons after GSO-projection. 

\subsection{Properties of the Lie superalgebra of physical states}

In this section we derive some properties of the Lie superalgebras $G$ and ${\tilde G}$. In the next section they will be used to show that $G$ is a generalized Kac-Moody superalgebra. 

We define a bilinear form on $G$ by setting
\[ \renewcommand{\arraystretch}{1.5}
\begin{array}{lcrcl} 
\la u, v \ra &=& (u,v)_H         & &        \mbox{if } u,v \in G_0 \\
\la u, v \ra &=& -(\tilde{u},v)_H & \mbox{with } 
                        u=X_{-1}\tilde{u}\, & \mbox{if } u,v \in G_1
\end{array} \]
and zero else. Note that by Proposition \ref{ncgp} this definition does not extend to $\tilde{G}$. We have
\begin{pp1}
The bilinear form $\la\, ,\,\ra$ pairs nondegeneratly $G(\al)$ with $G(-\al)$.
Furthermore $\la\, ,\,\ra$ is supersymmetric and invariant.
\end{pp1}
{\em Proof:}
The first part of the proposition is clear. It is also easy to see that 
$\la u, v \ra=\la v,u \ra$ for $u,v\in G_0$. Let $u,v \in G_1$. Then
\beann
\la u,v \ra &=&  -(c_{-2}\tilde{u}, X_{-1}\tilde{v}) \\
            &=&  -(c_{-2}\tilde{u}, X_{-1}b_1c_{-2}\tilde{v}) \\
            &=&  -(b_1X_{-1}c_{-2}\tilde{v}, c_{-2}\tilde{u}) \\
            &=&   (X_{-1}c_{-2}\tilde{v}, b_1c_{-2}\tilde{u}) \\
            &=&   (X_{-1}c_{-2}\tilde{v},\tilde{u}) \\
            &=&   (c_{-2}\tilde{v},X_{-1}\tilde{u}) \\
            &=&   (c_{-2}\tilde{v},u) \\
            &=&  - \la v,u \ra \, .
\eeann 
Hence $\la u, v \ra$ is supersymmetric.\\
To prove the invariance it is sufficient to consider two cases.
First let $u,v$ and $w$ be in $G_0$. Then
\beann 
\la [u,v],w \ra 
&=& - \la [v,u],w \ra \\
&=& - (-1)^{|v|} (c_{-2}b_1v_{-1}X_{-1}u,w) \\
&=&   (v_{-1}X_{-1}u , b_1c_{-2}w) \\
&=&   (v_{-1}X_{-1}u , w) \\
&=& - (X_{-1}u, v_{-1}^*w) \\
&=& - (X_{-1}u, v_{-1}w )
\eeann
We explain the last equality. Recall the definition
\[ v_{-1}^* = v_{-1}+\sum_{n\geq 1}\left(\frac{L_1^n}{n!}v\right)_{-n-1}\, . \]
For $n\geq 1$ we have
\[ L_1^nv = \{Q,b_2\}L_1^{n-1}v = Qb_2L_1^{n-1}v \]
Since $Q$ is a derivation this implies that only $v_{-1}w$ contributes to the product. Now
\beann 
\la u,[v,w] \ra &=& (-1)^{|v|}(c_{-2}u,b_1X_{-1}v_{-1}w) \\
                &=& -(b_1X_{-1}v_{-1}w,c_{-2}u) \\
                &=& -(X_{-1}v_{-1}w,u) \\
                &=& -(v_{-1}w,X_{-1}u) \\
                &=& -(X_{-1}u,v_{-1}w)
\eeann
proves the invariance in this case. Now let $u,v\in G_1$ and $w\in G_0$.
Then
\beann 
\la [u,v],w \ra &=& \la [v,u], w \ra \\
                &=& (-1)^{|v|}(c_{-2}b_1 v_{-1}u, w) \\
                &=& (v_{-1}u,b_1c_{-2}w) \\
                &=& (v_{-1}u,w)\\
                &=& (u,v_{-1}^*w) \\
                &=& (u,v_{-1}w)
\eeann
by the same argument as above and 
\beann
\la u,[v,w] \ra &=& - \la [v,w] ,u \ra \\
                &=& (-1)^{|v|}(c_{-2}b_1v_{-1}w,u ) \\
                &=& (v_{-1}w,b_1c_{-2}u) \\
                &=& (v_{-1}w,u) \\
                &=& (u,v_{-1}w) 
\eeann
showing the invariance in this case.
\hspace*{\fill} $\Box$\\

Simple calculation gives the bilinear form on the massless states
\begin{pp1}
Let $\al\in L^X$ with $\al^2=0$ and $\al\neq 0$. Then
\[ \renewcommand{\arraystretch}{1.3}
\begin{array}{rcl}
 \la \, \xi,\al\,|\, \zeta, -\al \, \ra &=& \xi_{\mu} \zeta_{\nu} g^{\mu \nu}\\
 \la \, X_{-1}\,|\,u,{-\textstyle \frac{3}{2}},\al \ra , 
                                |\,v,{-\textstyle \frac{1}{2}},-\al \, \ra  
                        &=&  u_{\bt}v_{\dot \gamma} C^{\bt \dot \gamma} 
\end{array} \]
For $\al =0$ we have 
\[  \la P^{\mu} , P^{\nu} \ra = g^{\mu \nu}\, . \]
\end{pp1}
Hence the map $j$ from $H=G(0)_0$ to ${\mathbb C}\otimes_{\mathbb Z}L^X$ defined by $P^{\mu}\mapsto x^{\mu}$ is an isometry. We have
\begin{pp1} 
Let $h\in H, \, x\in G(\al)$ and $y\in G(-\al)$. Then
\beann
{[} h,x {]} &=& (h,\al) x \\
{[} x,y {]} &=& \la x, y \ra j^{-1}(\al)
\eeann
\end{pp1}
{\em Proof:} From $b_0X_{-1}P^{\mu}=-x^{\mu}(-1)$ we get 
\be 
[P^{\mu}, x] = \{X_{-1}P^{\mu},x\}= -(b_0X_{-1}P^{\mu})_0 =\al^{\mu}x \, .
\ee
This proves the first statement. The second follows from the nondegeneracy of 
$\la\, ,\, \ra$. Let $h\in H$. Then
\[ \la\, [x,y]-\la x, y\ra j^{-1}(\al), h\ra=
   \la x,[y,h]\, \ra - \la x, y\ra (\al , j(h)) =0 \, .\]
\hspace*{\fill} $\Box$\\
We give a simple application of this result.
\begin{pp1}
$G$ is simple.
\end{pp1}
{\em Proof:} Let $I\neq 0$ an ideal in $G$. Then by Proposition 1.5 in 
\cite{K1} $I$ is graded by $L^X$. $I$ clearly has an element of degree $\al\neq 0$ in $L^X$. This implies that $I\cap H$ is nonempty. Let $h$ be a nonzero element in $I\cap H$. Since $\la \, , \, \ra$ is nondegenerate the elements $h'$  
with $\la h,h' \ra \neq 0$ span $H$. $I$ contains the root spaces $G(\al')$ corresponding to these $h'$. This implies that the $h'$ are in $I$. Hence $I\cap H=H$ and $I=G$. \hspace*{\fill} $\Box$\\

This implies that $G$ is the derived algebra of ${\tilde G}$.

The next proposition shows that the zero momentum states $P^{\mu}$ and 
$Q^{\dot{\al}}$ generate a representation of the $N=1$ supersymmetry algebra on $G$ and ${\tilde G}$. 
\begin{pp1}
The zero momentum states satisfy the following algebra.
\[ \renewcommand{\arraystretch}{1.5}
\begin{array}{c}
 [ P^{\mu} , P^{\nu} ] = [ P^{\mu},Q^{\dot{\al}} ] = 0 \\
 \{ Q^{\dot{\al}},Q^{\dot{\bt}} \} = 
         \frac{1}{\sqrt 2} (\Gamma_{\mu}C)^{\dot{\al}\dot{\bt}} P^{\mu}  
\end{array} \]
\end{pp1}
{\em Proof:} The first two products are clear. 
We do the calculations for the third product. 
\beann
\{ Q^{\dot{\al}},Q^{\dot{\bt}} \} 
&=& (b_0 S^{\dot{\al}}_{-1}c)_0S^{\dot{\bt}}_{-1}c \\
&=& - c_{-1} S^{\dot{\al}}_0 S^{\dot{\bt}} \\
&=& - c_{-1}\frac{1}{\sqrt 2} (\Gamma_{\mu}C)^{\dot{\al}\dot{\bt}} 
             (\psi^{\mu}_{-1}e^{-\phi}) \\
&=& \frac{1}{\sqrt 2} (\Gamma_{\mu}C)^{\dot{\al}\dot{\bt}} P^{\mu}
\eeann 
by Proposition \ref{susyv}.
\hspace*{\fill} $\Box$\\
Hence the spaces $G(\al)_0$ and $G(\al)_1$ are related by a supersymmetry transformation for $\al\neq 0$. This implies 
\begin{thp1}
Let $\al\neq 0$. Then the dimension of the root space $G(\al)_1$ is given by
\[ \mbox{dim }G(\al)_1 = c({\textstyle -\frac{1}{2}}\al^2) \, . \]  
\end{thp1}

The following technical result is proved by direct calculation.
\begin{pp1}
Let $\al$ and $\bt$ be proportional zero norm vectors in 
$L^X\backslash \{0\}$. 
Then their root spaces commute.           
\end{pp1}

Consider the real vector space spanned by the states 
\[    x^{\mu_1}(-1)_{-m_1}\ldots x^{\mu_i}(-1)_{-m_i}
      \psi^{\nu_1}_{-n_1}\ldots \psi^{\nu_j}_{-n_j} 
      e^{-\phi}_{-1}e^{\al}_{-1}c \]
where $j$ is odd and the $m$'s and $n$'s are positive integers and define the 
subspace $P(\al)_{-1}$ of states $v$ satisfying  
$L^M_0 v= \frac{1}{2}v$ and $L^M_nv=G^M_{n-\frac{1}{2}}v=0$ for $n\geq 1$.
Then $P(\al)_{-1}$ is $Q$-closed. We denote the image in $H(\al)_{-1,1}$ also 
by $P(\al)_{-1}$. Then the complexification of $P(\al)_{-1}$ is 
$H(\al)_{-1,1}$ (cf. e.g. \cite{P}).
Similarly the states 
\[    x^{\mu_1}(-1)_{-m_1}\ldots x^{\mu_i}(-1)_{-m_i}
      \psi^{\nu_1}_{-n_1-\frac{1}{2}}\ldots 
      \psi^{\nu_j}_{-n_j-\frac{1}{2}} 
      S^{\dot{\bt}}_{-1}e^{\al}_{-1}c        \]
where $j$ is even and the $m$'s and $n$'s are again positive integers
generate a real vector space and we define $P(\al)_{-\frac{1}{2}}$ as the 
subspace of states $v$ satisfying $L^M_0 v= \frac{5}{8}v$ and 
$L^M_nv=G^M_{n-1}v=0$ for $n\geq 1$. $P(\al)_{-\frac{1}{2}}$ is $Q$-closed and its image in $H(\al)_{-\frac{1}{2},1}$ generates this space.

Now choose a basis of spinors such that the $\Gamma$-matrices defined above 
satisfy the Majorana and the Weyl condition. Denote the sum of the 
$P(\al)_{-1}$ over $\al$ in $L^X$ as $P_{-1}$ and the sum of the 
$P(\al)_{-\frac{1}{2}}$ with $\al\neq 0$ as $P_{-\frac{1}{2}}$.
Then one can prove
\begin{pp1}
$P_{-1} \oplus P_{-\frac{1}{2}}$ is a real form of $G$. The bilinear form on $G$ induces a real bilinear form on this space. 
\end{pp1}

\subsection{The main theorems} \label{mrs}

In this section we state the main results of this paper.
\begin{thp1} \label{mt1}
The Lie superalgebra $G$ has the following properties.
\begin{enumerate}
\item[]
$G$ is graded by $L^X$. We define $H$ as the subspace of degree zero, so that
$G=H\oplus \bigoplus_{\al\neq 0}G(\al)$.
\item[]
$G$ has a invariant supersymmetric bilinear form that pairs nondegenerately
$G(\al)$ with $G(-\al)$.
\item[]
$H$ is an abelian selfcentralizing even subalgebra of $G$. $H$ is isometric to ${\mathbb C}\otimes_{\mathbb Z}L^X$ and $[h,x]=(h,\al)x$ for $h\in H$ and 
$x\in G(\al)$ so that the $G(\al)$ are eigenspaces of $H$. $H$ is called the Cartan subalgebra of $G$. 
\item[]
The root spaces $G(\al)=G(\al)_0+G(\al)_1$ are finite dimensional and
\linebreak  
$\mbox{dim}\, G(\al)_0 = \mbox{dim}\, G(\al)_1 = c(-\frac{1}{2}\al^2)$.
Hence the roots of $G$ are the nonzero $\al\in L^X$ with $\al^2\leq 0$.  
\item[]
$G$ is simple.
\item[]
$G$ is acted on by the $N=1$ superymmetry algebra.
\item[]
The root spaces of two proportional norm zero roots commute.
\item[]
$G$ is a generalized Kac-Moody superalgebra.
\end{enumerate}
\end{thp1}
\pagebreak[3]
{\em Proof:} Using the root space decomposition and the nondegeneracy of $\la \, ,\, \ra$ it is easy to see that $H$ is selfcentralizing. The last statement follows from Theorem \ref{ut}. As regular element one can take any nonzero negative norm vector. The roots of $G$ are all of infinite type. Note that in a Lorentzian lattice any two positive imaginary roots have inner product at most zero and zero only if they are both multiples of the same norm zero vector.
\hspace*{\fill} $\Box$\\

Now we consider the case that $L^X$ is the unique unimodular even Lorentzian 
lattice $II_{9,1}$. In this case we call $G$ the fake monster superalgebra.
The lattice $II_{9,1}$ can be described as the lattice of all points 
$(x_1,\ldots,x_9,x_{10})\in {\mathbb R}^{9,1}$ with all $x_i\in {\mathbb Z}$ or all $x_i\in {\mathbb Z}+\frac{1}{2}$ and which have integer inner product with $(\frac{1}{2},\ldots,\frac{1}{2},\frac{1}{2})$. In a Lorentzian lattice there are two cones of negative norm vectors. We define one of them as the positive cone. Then we have
\begin{thp1}
The simple roots of the fake monster superalgebra are the norm zero vectors in the closure of the positive cone of $II_{9,1}$. Their multiplicities as even and odd roots are equal to 8. 
\end{thp1}
{\em Proof:}
Let $K$ be the generalized Kac-Moody superalgebra with Cartan subalgebra 
${\mathbb C}\otimes_{\mathbb Z}II_{9,1}$ and simple roots as stated in the theorem. The simple roots of $G$ are determined by the root multiplicities because of the denominator formula. Hence it is sufficient to show that $G$ and $K$ have the same root multiplicities so that they are isomorphic. 
The denominator identity (cf. \cite{GN}) for $K$ is 
\be \prod_{\al\in\Delta_{+}}
   \frac{ (1-e(\al))^{ \mbox{mult}_0 \, \textstyle{\al} } }
        { (1+e(\al))^{ \mbox{mult}_1 \, \textstyle{\al} } }
        = e(\rho) \sum_{w\in W} \mbox{det}(w)w(T)  \ee
where $T$ is the sum 
\[ T=e(-\rho)\sum \epsilon(\mu)e(\mu) \]
 with $\mu$ running over all sums of imaginary simple roots and $\epsilon(\mu)=(-1)^n$ if $\mu$ is the sum of $n$ pairwise perpendicular imaginary simple roots, which are distinct unless they are odd and of norm zero, and $\epsilon(\mu)=0$ else. Note that $e(\al)$ and $e(\mu)$ are elements of the group algebra of $II_{9,1}$. We can calculate the right hand side of the identity, i.e. the denominator function of $K$. Since $K$ has no real simple roots the Weyl group of $K$ is trivial. The Weyl vector is zero. Hence we only have to calculate $\sum \epsilon(\mu)e(\mu)$. $\mu$ is the sum of pairwise perpendicular imaginary simple roots iff it is the multiple of a primitive norm zero vector in $II_{9,1}$ since two imaginary simple roots can only be are orthogonal if they are proportional. Taking care of the multiplicities of the simple roots as elements in the Cartan subalgebra we find
\[ \sum \epsilon(\mu)e(\mu) = 1 + \sum a(\lambda)e(\lambda) \]
where $a(\lambda)$ is the coefficient of $q^n$ of 
\[ \prod_{m\geq 0} \left( \frac{1-q^m}{1+q^m} \right)^8 =
   1-16q+112q^2-448q^3+\ldots \]
if $\lambda$ is $n$ times a primitive norm zero vector in the closure of the positive cone and zero else.   
Borcherds has shown (cf. \cite{B5}, Example 13.7) that  
\[ 1 + \sum a(\lambda)e(\lambda) = \prod_{\al\in\Delta_{+}}
   \frac{ (1-e(\al))^{ c(-\frac{1}{2}\al^2) }}
        { (1+e(\al))^{ c(-\frac{1}{2}\al^2) }} \]
This implies that the multiplicity of a root $\al$ of $K$ is given by $c(-\frac{1}{2}\al^2)$. \hspace*{\fill} $\Box$\\ 

\begin{cp1}
The denominator identity of the fake monster superalgebra is 
\[ \prod_{\al\in\Delta_{+}}
   \frac{ (1-e(\al))^{ c(-\frac{1}{2}\al^2) }}
        { (1+e(\al))^{ c(-\frac{1}{2}\al^2) }} =
    1 + \sum a(\lambda)e(\lambda) \]
where $a(\lambda)$ is the coefficient of $q^n$ of 
\[ \prod_{m\geq 0} \left( \frac{1-q^m}{1+q^m} \right)^8 =
   1-16q+112q^2-448q^3+\ldots \]
if $\lambda$ is $n$ times a primitive norm zero vector in the closure of the positive cone and zero else. 
\end{cp1}

The denominator function of the fake monster superalgebra is the right hand side of the denominator identity. It is an automorphic form of weight $4$ for a subgroup of $O_{10,2}({\mathbb R})$ (cf. \cite{B5}, Example 13.7). 

Finally we describe the construction of the fake monster superalgebra $G$ by generators and relations. Let $\{ \al_i \,|\, i\in I\}$ the set of simple roots of $G$. Taking $8$ even and $8$ odd copies of each $\al_i$ we can calculate the Cartan matrix of $G$. Then we define as described in section 3 the superalgebra $\hat{G}$ by generators and relations.
The subalgebra $\hat{H}$ of $\hat{G}$ is the direct sum of a $16^2=256$ dimensional space for each simple root $\al_i$. 
The center $C$ of $\hat{G}$ has index $10$ in $\hat{H}$. The fake monster superalgebra is the quotient of $\hat{G}$ by $C$.

\section*{Acknowledgments}

I thank R. E. Borcherds for stimulating discussions. Furthermore I want to thank N. Berkovits and the referee for useful comments.

\end{document}